\documentclass{article}

\usepackage[pctex32]{graphics}
\usepackage{amssymb}
\usepackage{latexsym}
\usepackage{epsfig} 
\usepackage{multirow} 
\usepackage{pstricks}
\usepackage{graphicx}
\usepackage{float}
\usepackage{epstopdf}
\usepackage{mathrsfs}
\usepackage{subcaption}
\usepackage{soul}

\usepackage{amsmath}
\usepackage{amssymb}
\usepackage{relsize}
\usepackage{multirow}
\usepackage{color}
\usepackage{graphicx}
\usepackage{multicol} 

\usepackage{times}
\usepackage{epstopdf}
\newcommand{\comm}[1]{}
\usepackage{fancyhdr}
\usepackage{setspace}
\usepackage{pstricks} 
\newcommand{\R}{{\mathbb R}}

\newcommand{\mW}{{\mathsf W}}

\newcommand{\mP}{{\mathsf P}}

\newcommand{\mT}{{\mathsf T}}
\newcommand{\mC}{{\mathsf C}}
\newcommand{\mI}{{\mathsf I}}
\newcommand{\mA}{{\mathsf A}}
\newcommand{\mV}{{\mathsf V}}
\newcommand{\mE}{{\mathsf E}}

\newcommand{\mL}{{\mathsf L}}
\newcommand{\mQ}{{\mathsf Q}}

\newcommand{\mD}{{\mathsf D}}

\newcommand{\mU}{{\mathsf U}}
\newcommand{\mG}{{\mathsf G}}

\newcommand{\mH}{{\mathsf H}}
\newcommand{\mSigma}{{\mathsf \Sigma}}

\newcommand{\mLambda}{{\mathsf \Lambda}}

\newcommand{\vertiii}[1]{{\left\vert\kern-0.25ex\left\vert\kern-0.25ex\left\vert #1 
    \right\vert\kern-0.25ex\right\vert\kern-0.25ex\right\vert}}

\title{Dictionary learning methods for brain activity mapping with MEG data}
\author{Daniela Calvetti \and Erkki Somersalo}
\date{Case Western Reserve University\\
Department of Mathematics, Applied Mathematics, and Statistics\\
Cleveland, OH, USA}

\begin{document}

\maketitle

\begin{abstract}

A central goal in many brain studies is the identification of those brain regions that are activated during an observation window that may correspond to a motor task, a stimulus, or simply a resting state. While functional MRI is currently the most commonly employed modality for such task, methods based on the electromagnetic activity of the brain are  valuable alternatives because of their excellent time resolution and of the fact that the measured signals are directly related to brain activation and not to a secondary effect such as the hemodynamic response. In this work we focus on the MEG modality, investigating the performance of a recently proposed Bayesian dictionary learning (BDL) algorithm for brain region identification. The partitioning of the source space into the 148 regions of interest (ROI) corresponding to parcellation of the Destrieux atlas provides a natural determination of the subdictionaries necessary for the BDL algorithm. We design a simulation protocol where a small randomly selected patch in each ROI is activated, the MEG signal is computed and the inverse problem of active brain region identification is solved using the BDL algorithm. The BDL algorithm consists of two phases, the first one comprising dictionary compression and Bayesian compression error analysis, and the second one performing  dictionary coding with a deflated dictionary built on the output of the first phase, both steps relying on Bayesian sparsity promoting computations. For assessing the performance, we give a probabilistic interpretation of the confusion matrix, and consider different impurity measures for a multi-class classifier.
\end{abstract}

\section{Introduction}

Understanding the connectivity and the functional organization of the human brain is an important aspect of brain research, including studies on  consciousness, sleep,  autism, different neuro-degenerative diseases, language development,  intelligence and more. Understanding the network structure of the brain may also help understanding how epileptic seizures are triggered.
An essential part of this endeavor is parcellation of the brain into regions that can be identified either based on the brain anatomy or by the function. While the overall functional role of some brain areas, for example motor cortex, olfactory bulb, or visual cortex, is well understood, a finer parcellation and
the association of detected activation with a specific function
is challenging because the lack of definite ground truth, possibly leaving room for different interpretations. In summary, although automatic anatomy-based parcellation techniques have been proposed in the literature, see, e.g., \cite{lancaster2000automated,desikan2006automated,destrieux2010automatic},  unambiguous identification of functional areas is more challenging, as pointed out, e.g.,  in \cite{salehi2020there}. In particular, the functional organization of the brain may be hierarchical and multi-scale in nature as some studies indicate \cite{doucet2011brain,hilgetag2020hierarchy}.

Most functional brain mapping studies are based on PET or fMRI imaging, see, e.g. \cite{raichle2009brief} for a brief history of functional brain mapping. Currently, fMRI is the main imaging tool in cognitive neuroscience, however, its limitations are known and have been discussed in the literature \cite{logothetis2008we}. In particular, it is well understood that the fMRI modality is sensitive to the oxygenation level of the blood, which does not directly correspond to cerebral metabolic rate (CMR) of oxygen, raising the question of its direct correspondence to brain activity. The interpretation of the BOLD data as a direct proxy for neuronal firing activity is further complicated by the fact that inhibitory activity increases CMR but is not associated with neuronal firing, implying that fMRI may be unable to distinguish excitation and inhibition unequivocally. Magnetoencephalography (MEG), on the other hand,  seeks to estimate directly the neuronal electrophysiological synaptic activity, thereby offering a viable alternative for cognitive neuroimaging. One question that needs to be addressed with MEG is whether it can correctly identify active brain regions, which is the main topic of this study. 

Unlike in the standard MEG inverse problem of source localization \cite{baillet2001electromagnetic,hamalainen1993magnetoencephalography},
we consider the problem of identification of active brain regions based on MEG data, i.e., singling out a few brain regions containing the electric currents inducing the magnetic field withing the observation accuracy. In particular, different sources in the same brain region are considered equivalent, and the goal is to identify as few as possible regions so that the data are satisfactorily explained. The proposed algorithm therefore involves sparse coding in terms of the brain regions specified by a given brain atlas. In the context of MEG and EEG, sparsity-promoting methods have been considered for decades, either based on appropriately chosen penalties, see, e.g., \cite{uutela1999visualization,krishnaswamy2017sparsity,pascarella2019inversion}, or more in line with the current analysis, using Bayesian sparsity-promoting priors \cite{calvetti2015hierarchical,calvetti2019hierachical,calvetti2023ias}. In the current setting, rather than looking for a sparse representation of the activated source in terms of single dipoles, we seek to explain the data in terms of a few groups of sources, yielding a group sparsity problem.  Group sparsity-promoting penalty methods have been previously proposed in the context of MEG, see, e.g., \cite{lim2017sparse}. In this work, the groups sparsity methods 
are developed further in the Bayesian framework, and in particular, in the context of dictionary learning, where the dictionary interpretation is guided by an underlying brain parcellation. One way to lessen the computational burden is to first perform a compression of the dictionary, which can be regarded as a model reduction technique and, following the ideas of Bayesian modeling error analysis \cite{kaipio2006statistical,kaipio2007statistical}, augment the compressed dictionary with the associated compression error. A statistically motivated group sparsity prior, similar to the anatomical prior introduced in \cite{calvetti2015hierarchical}, is introduced to improve the success rate of the algorithm. Following  \cite{bocchinfuso2024bayesian}, the compressed dictionary matching algorithm is used to deflate the original dictionary retaining only the portion that may be pertinent to explain the MEG data. This leads to a two-phase algorithm based on recent ideas of Bayesian sparsity promotion. The proposed algorithm is utilized in a large set of numerical simulations designed to test how well MEG can be used to detect activation in the different brain regions.

\section{Brain activity mapping}

A standard model for the electromagnetic brain activity represents the impressed current associated with the post-synaptic ion currents in the brain as a discrete sum of point-like dipoles,
\[
 J = \sum_{k=1}^N q_k,
\]
where $q_k\in\R^3$ is the moment of the $k$th dipole located at $p_k\in\R^3$. The ensemble of dipole positions constites the source space $S = \{p_k:1\leq k\leq N\}$. In the present article we consider the magnetoencephalography (MEG) data, comprising the noisy observation of the magnetic flux density components by a number of magnetometers, and possibly the gradiometer data, i.e., gradients of the flux densities in the direction perpendicular to the magnetometer axis. Let $m$ denote the number of data channels. The quasi-static Maxwell's equations allow us to represent the measured components of the magnetic flux density and/or their tangential gradients as
\begin{equation}\label{forward}
 b = \sum_{k=1}^N \mL_k q_k,
\end{equation}
where $b\in\R^m$ is the vector of all observed scalar quantities, and $\mL_k\in \R^{m\times 3}$ is the lead field matrix accounting for the geometry and conductivity distribution of the head; for details of modeling and computations, see, e.g., \cite{ilmoniemi2019brain,wolters2004efficient}.

While the exact location of the brain activity is crucial in applications such as planning epilepsy surgery, to understand the time course of brain functions during a task or in the resting state, or to map the functional connectivity of the brain, it is sufficient to identify the active brain region within a given brain atlas. In that context, it is natural to approach the problem from the point of view of dictionary learning and dictionary matching.

\subsection{Dictionary matching with group sparsity}

For simplicity, consider the forward model (\ref{forward}) under the simplifying hypothesis that the dipole directions are fixed to the preferential direction, typically normal to the cortex. Denoting by $\widehat q_k\in{\mathbb S}^2$ the unit vector determining the preferential direction at the source space location $p_k$, we write the forward model as
\[
  b = \sum_{k=1}^N \big(\pm|q_k|\big) \big(\mL_k \widehat q_k\big).
\]
At this point, the effects of observation noise and
unspecified modeling uncertainties are not included in the representation.
Furthermore, in the present context, the dipole and signal strengths are of secondary interest. Therefore, given an observation vector $b$, we start by considering the {\em dictionary matching problem} of finding a representation 
\begin{equation}\label{dictionary matching}
 y = \sum_{k=1}^n d_k x_k = \mD x,
\end{equation}
where 
\[
 y = \frac b{\|b\|}, \quad d_k = \frac{\mL_k \widehat q_k}{\|\mL_k \widehat q_k\|}, 
\]
ensuring that the columns of the matrix $\mD$ are scaled to have unit Euclidean norm and the coefficients $x_k$ are proportional to the dipole amplitudes. The recasting (\ref{dictionary matching}) of the source identification problem amounts to matching the query $y$, the scaled MEG data, with the dictionary atoms $d_k$ encoding the dipole contributions. 

A physiological atlas of the brain comprising $L$ regions, given in terms of disjoint partitioning of the index set
\[
 I =\big\{1,2,\ldots,n\big\} = \bigcup_{\ell =1}^L I_\ell, \quad I_\ell\cap I_{\ell'} = \emptyset \mbox{ for $\ell\neq \ell'$,}
\]
can be used to obtain a natural parcelling of the dictionary atoms so that we can reformulate the dictionary matching problem as 
\[
 y = \sum_{\ell=1}^L\left(\sum_{k\in I_\ell} d_k x_k\right) =\sum_{\ell  = 1}^L \mD^{(\ell)} x^{(\ell)}.
\]
The columns $d_k$, $k\in I_\ell$ of the matrix $\mD^{(\ell)}\in\R^{m\times n_\ell}$ are the atoms in the $\ell$th subdictionary, and  the vector $x^{(\ell)} \in\R^{n_\ell}$ contains the corresponding weights. If the brain activity at a given time instance is concentrated in only few brain regions, it is reasonable to seek a  solution such that the vector $\big(\|x^{(1)}\|, \ldots, \|x^{(L)}\|\big)$ is sparse. This is tantamount to a {\em group sparsity} prior belief. In the literature, group sparsity has often been promoted through group LASSO regularization \cite{yuan2006model}, leading to the  minimization problem 
\[
 \min_x \|y - \mD x\|^2 + \lambda \sum\|x^{(\ell)}\|,
\]
where $\lambda>0$ is a properly chosen Lagrange multiplier that can be interpreted as a Tikhonov type regularization parameter. Our approach to the promotion of group sparsity is through computationally efficient hierarchical models. More specifically, we propose a two-phase classifier algorithm based on a Bayesian model for dictionary compression and deflation. In the following section, we outline the structure of the algorithm. The details of its implementations will be presented in section 3.

\subsection{A two-phase classifier: Dictionary compression and deflation}

In dictionary matching and learning algorithms it is common to compress the subdictionaries to summarize their content in terms of a few feature vectors thus reducing the degrees of freedom. Dictionary compression often is based on a low-rank approximation of each subdictionary 
\[ \mD^{(\ell)} \approx \mW^{(\ell)} \mH^{(\ell)}, \; \mW^{(\ell)} \in \R^{n \times k}, \; \mH^{(\ell)} \in \R^{k \times p_\ell} \]  using, e.g., a truncated singular value decomposition (SVD) or non-negative matrix factorization (NMF) \cite{aharon2006k,gillis2020nonnegative}.
The $\ell$th reduced dictionary $\mW^{(\ell)}$ comprises of few feature vector columns that approximately span the range of the full subdictionary $\mD^{(\ell)}$. The goal of the sparse compressed dictionary learning is to find a representation
\begin{equation}\label{obs model}
 y = \sum_{\ell = 1}^L \mW^{(\ell)} z^{(\ell)} +e
\end{equation}
in terms of as few subdictionaries as possible such that the model discrepancy $e$ is small. 

In the formulation of the dictionary matching problem in a Bayesian framework, $y$, $z^{(\ell)}$ and $e$
in (\ref{obs model}) are modeled as realizations of the random variables $Y$, $Z^{(\ell)}$ and $E$, respectively, with $E$ representing the additive uncertainty arising from the measurement error of $Y$ and the dictionary compression error (DCE) due to the passage from the full dictionary to the compressed dictionary 
\[
\mW = [\mW^{(1)}, \ldots, \mW^{(L)}].
\]
We refer to $E$ simply as DCE and model it a Gaussian random variable,
\[
 E \sim{\mathcal N}(\mu, \mC);
\]
where the estimation of the mean and the covariance will be discussed in the following section. Likewise, we use a conditionally Gaussian prior model for each $Z^{(\ell)}$. Introducing scalar random variables $\Theta_\ell$, we write
\[
 Z^{(\ell)}\mid (\Theta_\ell = \theta_\ell) \sim {\mathcal N}(0,\theta_\ell \mG^{(\ell)}),
\]
where $\mG^{(\ell)}$ is a symmetric positive definite matrix to be specified later in the next section.
Finally, we observe that the variance parameters are crucial for correctly identifying the activity regions. Indeed, the conditional prior density of the random variable $Z = \big(Z^{(1)}, \ldots,Z^{(L)}\big)$ conditioned on $\Theta = \theta\in\R^L$ is
\[
 \pi_{Z\mid\Theta}(z\mid\theta) \propto \prod_{\ell = 1}^L \theta_\ell^{-r_\ell/2} {\rm exp}\left( - \frac 12 \frac{\|z^{(\ell)}\|_{(\ell)}^2}{\theta_\ell}\right)
 ={\rm exp}\left( - \frac 12 \sum_{\ell=1}^L \frac{\|z^{(\ell)}\|_{(\ell)}^2}{\theta_\ell}
 -\frac 12 \sum_{\ell=1}^L r_\ell \log\theta_\ell\right),
\]
where
\[
 \|z^{(j)}\|_{(j)}^2 = \big(z^{(j)}\big)^\mT \big(\mG^{(j)}\big)^{-1} z^{(j)}.
\]
To explain the data in terms of as few subdictionaries as possible, it is desirable a priori that the vector 
\[
 \vertiii{z} = \big(\|z^{(1)}\|_{(1)|},
 \|z^{(2)}\|_{(2)|},\ldots,\|z^{(L)}\|_{(L)}\big)
\]
is sparse, or compressible, meaning that most of the components are below some threshold value or, alternatively, most of the variances $\theta_\ell$ are small and only a few are of large magnitude. We do not know a priori which entries are significantly different from zero as this information is tantamount to knowing which brain regions are active. Rather, the significant variances will be determined a posteriori by the data. To that end, we introduce a hyperprior model for the random variable $\Theta$ expressing our belief that the components are mutually independent and identically distributed, and obeying a fat tailed distribution law that favor small values with occasional large outliers, for example a Generalized gamma distribution model,
\[
 \Theta_j \sim\Theta^{r\beta_j-1}{\rm exp}\left( - \left(\frac{\theta_j}{\vartheta_j}\right)^r\right).
\]
As explained in detail in the following section, we set up a Bayesian posterior model and seek for the Maximum A Posteriori (MAP) estimate for the pair $(Z,\Theta)$,
\begin{equation}\label{MAP}
(z^*,\theta^*) = \underset{(z,\theta)}{\rm argmax}\, \pi_{Z,\Theta\mid Y}(z,\theta\mid y).
\end{equation}
Finally, to identify a single dominant active brain region, we select the subdictionary with maximal variance,
\[
 i^* = \underset{1\leq i\leq L}{\rm argmax}\, \theta_i^*.
\]
The dictionary learning algorithm described above leads reliably to a compressible solution, most of the variance parameters $\theta_\ell$ being negligible. However, when used to find a single brain region as described above, the winner-takes-all strategy will fail even if its variance $\theta_\ell$ is relatively large, thus indicating a plausible candidate, but not the largest. As shown in  \cite{bocchinfuso2024bayesian}, the success of classification algorithm described above can be improved significantly by replacing the winner-takes-all algorithm by a dictionary deflation step followed by a sparse dictionary coding step with a deflated dictionary.

More precisely, once the MAP estimate (\ref{MAP}) has been computed, we identify the subset of  indices
\[
 I_{\rm defl}(p) = \{\ell \mid \theta_\ell > p\, \theta_{i^*}\}, \;0<p<1, 
\] 
where $p$ is a threshold parameter. The set $I_{\rm defl}(p)$ points to the few brain regions that the compression algorithm has identified as relevant for explaining the data, and is used to define the deflated dictionary
\begin{equation}\label{deflated}
 \mD^{\rm defl} = \left[\begin{array}{ccc} \mD^{(\ell_1)} & \cdots & 
 \mD^{(\ell_q)}\end{array}\right], \quad I_{\rm defl}(p) =\{\ell_1,\ldots,\ell_q\}.
\end{equation}
The second step in the classification algorithm is to find an optimal representation of the data vector in terms of a few atoms of the deflated dictionary,
\[
 y = \mD^{\rm defl} x + e', \quad \mbox{$x$ sparse,}
\]
where the vector $e'$ represents the observation uncertainty in the data, by minimizing the discrepancy while promoting sparsity for the vector $x$. The details of solution of the deflated dictionaty learning problem can be found in the next section. Having computed the coefficient vector $x$, we apply again a winner-takes-all classification strategy, associating the data $y$ with the subdictionary in $I_{\rm red}(p)$ that contributes most to its explanation. The criterion for the best match is done similarly as in the first step, by identifying the component with maximal variance in a Bayesian hierarchical model.

The strength of the two-phase algorithm is that, on the one hand, the first phase is more robust and less sensitive to noise than a sparsity-based dictionary matching with the full dictionary, as the feature vectors summarizing the subdirectories  typically are not a good match for noise. On the other hand, the lack of fine details in the compressed dictionary is compensated for in the second step,  where the dictionary matching problem is solved again with a deflated dictionary obtained by removing from the original dictionary those atoms deemed irrelevant in the first phase. The computational efficiency of this algorithm stems from the fact that both phases work with significantly reduced dictionaries.

\subsection{Performance assessment quantification}\label{sec:assessing}

Typically the performance of dictionary classifiers are tested by on a data sets consisting of inputs with known class label, ${\mathscr T}_{\rm test} = \big\{(y^{(k)},c_k)\big\}_{k=1}^K$. In the brain activity identification from MEG data, the testing will be done on simulated data, because of the general unavailability of  ground truth with experimental data.
Here, we limit the discussion to simulated activity in a single brain region, detailed being discussed with the computed examples.

Let $\mC$ denote the $L\times L$ confusion matrix, where the entry $\mC_{ij}$ contains the number of times when class $j$ has been assigned to data arising from source in class $i$.
A row-wise interpretation of  the confusion matrix gives the frequency at which each active brain region is assigned to data with a given label. Particularly interesting from the point of view of the inverse problem of region identification is the column-wise interpretation. The entries of the $j$th column of $\mC$ scaled by its 1-norm,
\[
\mP_{ij} = \frac{\mC_{ij}}{\sum_{i'=1}^L \mC_{i'j},} \quad \mbox{$j$ given,}
\]
indicate the frequencies of the true activity regions $i$, given that the algorithm has identified the $j$th region as the winner. 
Therefore we can interpret the entries of the column-normalized confusion matrix as conditional probabilities,
\[
 \mP_{ij} = P(\mbox{region $i$ is active}\mid \mbox{region $j$ was identified}), \quad 1\leq j\leq N.
\]
Similarly, the entries of the row-normalized confusion matrix,
\[
\mQ_{ij} = \frac{\mC_{ij}}{\sum_{j'=1}^L \mC_{i j'},}\quad \mbox{$i$ given,}
\]
admit the conditional probability interpretation,
\[
 \mQ_{ij} = P(\mbox{region $j$ is identified}\mid \mbox{region $i$ is active}), \quad 1\leq j\leq N.
\]
The matrix $\mP$ can be thought of as a generalization of {\em recall}
defined for binary classifiers, and the matrix $\mQ$ as a generalization of the {\em precision}.

To assess separately the performance of each phase, we denote by $\mC^1$ the confusion matrix corresponding to the compressed dictionary step with the winner-takes-all classification, referred to as Phase I, and by $\mC^2$ the confusion matrix after the completion of the deflated dictionary step, referred to as Phase II, and we use the notations $\mP^1$ and $\mP^2$ for the respective column-normalized versions.

Borrowing from the data science methodology, the performance of the algorithm associated to a specific brain region can be quantified by different impurity indices.  Given the matrix $\mP$, the misclassification rate $MCR_j$ gives the probability that the proposed classification to the $j$th class is incorrect,
\[
 MCR_j = \sum_{i\neq j}^L P_{ij} = 1 - P_{jj}.
\]
While the $MCR$ is a good measure for the success rate, it is insensitive to the spread of the  erroneous classifications. More informative measures of the spread are the Gini index $G_j$ and the cross entropy $S_j$,
\begin{equation}\label{impurity}
G_j =\sum_{i=1}^LP_{ij}(1-P_{ij}),\quad S_j =-\sum_{i=1}^L P_{ij}\log P_{ij}.
\end{equation}

We will use these impurity measures to assess the performance of the algorithm, as well as to assess to what extent the second phase improves/deteriorates the classification results.

\section{Algorithmic details}

This section gives the details of the implementation of the two-phase algorithm
that follows the lines of the classifier in \cite{bocchinfuso2024bayesian}. 

\subsection{Dictionary compression}\label{sec:compression}

We review here a dictionary compression technique that used ideas from the Bayesian inverse problems paradigm. We start by writing the singular value decomposition (SVD)  of each subdictionary matrix,
\[
 \mD^{(\ell)} = \mU^{(\ell)} \mSigma^{(\ell)} \left(\mV^{(\ell)}\right)^\mT,
\]
where $\mU^{(\ell)} \in\R^{m\times m}$ and $\mV^{(\ell)}\in \R^{n_\ell\times n_\ell}$ are orthogonal matrices and $\mSigma\in\R^{m\times n_\ell}$ is diagonal with non-negative diagonal entries $\sigma_j^{(\ell)}$ appearing in decreasing order. For each $\ell$, we then write a low-rank approximation of the subdictionary by discarding singular values that are below a given threshold value $\tau$,  $0<\tau<1$, relative to the largest one, 
\[
 \mD^{(\ell)} \approx \sum_{j=1}^{r_\ell} \sigma_j^{(\ell)} u^{(\ell)}_j \big(v^{(\ell)}_j\big)^\mT = \mW^{(\ell)} \mH^{(\ell)},\quad \sigma^{(\ell)}_{r_\ell}\geq\tau\sigma_1^{(\ell)}>\sigma^{(\ell)}_{r_\ell + 1},
\]
where the low rank factors are
\[
 \mW^{\ell)} = \left[\begin{array}{ccc} u^{(\ell)}_1 & \cdots & u^{(\ell)}_{r_\ell}\end{array}\right], \quad
 \mH^{(\ell)} = \left[\begin{array}{c} \sigma_1^{\ell)}\big(v^{(\ell)}_1\big)^\mT \\ \vdots \\
 \sigma_{r_\ell}^{(\ell)}\big(v^{(\ell)}_{r_\ell}\big)^\mT
\end{array}\right].
\]
We refer to the columns of $\mW_\ell$ as the feature vectors, or the code book, of the $\ell$th subdictionary, and the low-rank approximation is effectively a subdictionary compression. We write
\[
 \mD^{(\ell)} = \mW^{(\ell)} \mH^{(\ell)} + \mE^{(\ell)},
\]
where the matrix $\mE^{(\ell)}$ accounts for 
the compression error. We point out that while the low rank approximation was based on the SVD, the considerations about dictionary compression hold regardless of the low rank approximation method used. Repeating the procedure for all subdictionaries, we have
\begin{eqnarray*}
\mD &=&  \left[\begin{array}{cccc}
		\mA ^{(1)}	& \mD^{(2)} &\cdots &\mA ^{(L)}
	\end{array}\right] \\  &=&		 
\left[ \begin{array}{cccc} 
		\mW ^{(1)} \mH^{(1)}  & \mW^{(2)}\mH^{(2)} &\cdots 	&\mW ^{(L)} \mH^{(L)}
				\end{array} \right] + 
				\left[\begin{array}{cccc} \mE ^{(1)} & \mE^{(2)} &\cdots &\mE ^{(L)}
	\end{array}\right] \\
		&=&	\mW \mH + \mE,
\end{eqnarray*}
where
\[ \mW = 
	\begin{bmatrix}
		\begin{array}{ccc}
			\mW ^{(1)}	&\dots 	&\mW ^{(L)}
		\end{array} 
	\end{bmatrix}, \quad \mH =  \begin{bmatrix}			
		\mH^{(1)}	&0 	&\dots 	&0 \\
		0		&\mH^{(2)}	&0	&\vdots \\
		\vdots	&0	&\ddots	&0 \\
		0	&\dots	&0 	 	&\mH^{(L)}
	\end{bmatrix}. 
\]
The matrix $\mW$ is a compression of the  dictionary, and the goal of the first phase is to approximate a given input $y$ in terms of the atoms of $\mW$, which are the feature vectors of all the subdictionaries. 

We formulate the compressed dictionary matching problem in the Bayesian framework by restating it as
\[
 Y = \mW Z + E + S,
\]
where is a $Y$ random variable representing the input, $Z$ is a random variable of the dictionary coefficients, $E$ represents the dictionary compression error DCE, and $S$ accounts for the noise in the data. This model is the basis of the likelihood model of the problem. 
To analyze the DCE, consider an arbitrary  dictionary atom $d_k$, the $k$th column of $\mD$, and find its best approximation $\mW z_k$ in terms of the compressed dictionary. 
We decompose the atom as 
\begin{eqnarray*}
d_k  & = & \mW z_k + (d_k - \mW z_k )\\
& = & \mW z_k + e_k, \quad 1 \leq k \leq N.
\end{eqnarray*} 
If $d_k$ contains no measurement error, the discrepancy
\[
 e_k = d_k - \mW z_k , \quad 1\leq k\leq N,
\]
represents only the dictionary compression error. To encode the compression error in terms of a probability density, we assume that the vectors $e_k$ are independent realizations from the density, and use a Gaussian approximation of the form
\[
E  \sim \mathcal{N}(\mu_{\rm DCE}, \mC_{\rm DCE}).
\]
where the mean and covariance are estimated from the sample,
  \begin{eqnarray*}
 \mu_{\rm DCE} &=& \frac{1}{N} \sum e_k, \\
\mC_{\rm DCE} &=& \frac{1}{N-1} \sum (e_k - \mu_{\rm DCE}) (e_k - \mu_{\rm DCE})^\mT. 
\end{eqnarray*}
Assuming that the measurement error $S$ is independent of the compression error, and modeling it as scaled white noise,
\[
 S\sim{\mathcal N}(0, \delta^2 \mI_m),
\]
we can write the likelihood model for $Y$ as
\[
 \pi_{Y\mid Z}(y\mid z) \propto {\rm exp}\left( -\frac 12 \|\mC^{-1/2}(\mW z + \mu_{DCE} - y)\|^2 \right),
\]
where
\[
 \mC = \mC_{\rm DCE} + \delta^2 \mI_m.
\] 

Consider now the prior model.
To set up the model for $Z$, we proceed hierarchically. We begin by partitioning the random variable $Z$ in terms of the subdictionary components
\[
 Z = \left[\begin{array}{c} Z^{(1)} \\ \vdots \\ Z^{(L)} \end{array}\right],
\quad Z^{(k)} \in\R^{r_k},
\]
and introducing for each $Z^{(k)}$ a conditionally Gaussian structural prior density of the form
\[
 \pi_{Z^{(k)}\mid\Theta_k}(z^{(k)}\mid\theta_j) \propto
 \theta_j^{-r_k/2} {\rm exp}\left(-\frac 12 \frac{\| z^{(k)}\|_{(k)}^2}{\theta_j}\right), \quad
 \|z^{(k)}\|^2_{(k)}  = \left(z^{(k)}\right)^\mT
\left(\mG^{(k)}\right)^{-1}z^{(k)}.
\]
The covariance matrix $\mG^{(k)}$ is constructed as follows. First, consider the matrix $\mH^{(k)} \in \R^{r_k\times n_k}$, where $n_k$ is the number of atoms in the original dictionary that belong to the $k$th subdictionary. The columns of $\mH^{(k)}$ can be considered to represent typical coefficient vectors for matching data from the $k$th subdictionary, however, since we do not know a priori if the input vector $Z$ is in the $k$th subdictionary, we want to extract from $\mH^{(k)}$ only information about the direction of its columns. To that end, we compute the SVD of the matrix $\mH^{(k)}$, 
\[
 \mH^{(k)}= \mQ^{(k)} \mLambda^{(k)} \big[\mW^{(k)}\big]^\mT,  \quad \mLambda^{(k)} = {\rm diag}\big([\lambda^{(k)}_1,\ldots,\lambda^{(k)}_{p_k},0,\ldots,0]\big),
\] 
 and define
 \begin{eqnarray*}
\mG^{(k)}&=&\left(\frac{1}{\lambda_1^{(k)}}\right)^2 \mQ^{(j)} \mLambda^{(k)} \big[\mLambda^{(k)}\big]^\mT \big[\mQ^{(k)}\big]^\mT + \epsilon\, \mI \\
&=&  q_1^{(k)} \big[q_1^{(k)}\big]^\mT + \sum_{\ell=2}^{p_k} \left(\frac{\lambda_\ell^{(k)}}{\lambda_1^{(k)}}\right)^2 q_\ell^{(k)} \big[q_\ell^{(k)}\big]^\mT + \varepsilon\, \mI_{r_k},
\end{eqnarray*}
where $\varepsilon>0$ is a small parameter and $\varepsilon\, \mI_{r_k}$ is a regularization term added to guarantee positive definiteness. Intuitively, the first singular vector $q_1^{(k)}$ points in the direction in $\R^{r_k}$ of the cloud defined by the columns of $\mH^{(k)}$, while the subsequent singular vectors 
$q_\ell^{(k)}, \; \ell >1$ point in the directions where the point cloud has a spread proportional to $\big(\lambda_\ell^{(k)}\big)^2$. In this way, the prior expresses the belief that the direction of the coefficient vector $Z^{(k)}$ is expected to be guided by the cone along $q_1^{(k)}$ and  with the opening determined by the spread of the point cloud in the singular directions, while its amplitude is determined by the hyperparameter $\theta_k$.

We now proceed to define the prior distribution of the random variable modeling the parameter vector $\Theta\in\R^K$. Following the ideas of promoting sparsity via hierarchical models, we model the components $\Theta_j$ as mutually independent random variables with a heavy tailed distribution over the positive real axis, for example a generalized gamma distribution,
\[
 \Theta_k \sim\Theta^{r\beta_k-1}{\rm exp}\left( - \left(\frac{\theta_k}{\vartheta_k}\right)^r\right), \quad \theta_k>0,
\]
with $r\neq 0$, $\beta_k>0$ and $\vartheta_k>0$. The significance of $\beta_k$ and $\vartheta_k$ and will be discussed in the context of the numerical simulations. Heavy tailed distributions are a natural choice for the variances that we expect to be mostly small, with a possibility for large outliers. Motivated by these considerations, we set the prior model for the pair $(Z,\Theta)$ to be
\begin{eqnarray*}
 \pi_{Z,\Theta}(z,\theta) &=& \pi_{Z\mid\Theta}(z\mid\theta)\pi_\Theta(\theta) \\
 &\propto& {\rm exp}\left(-\sum_{k=1}^L\frac{\|z^{(k)}\|^2}{2\theta_j} -\sum_{k=1}^L\left(\frac{\theta_k}{\vartheta_k}\right)^r - \sum_{k=1}^L\left(\frac{r_k}{2} +1 -r\beta_k \right)\log\theta_k\right).
\end{eqnarray*}
Combining likelihood and prior according  to  Bayes' formula yields a posterior model for the pair $(Z,\Theta)$ of the form
\begin{eqnarray*}
    \pi_{Z,\Theta\mid Y}(z,\theta\mid y) &\propto&  \pi_{Y\mid Z}(y\mid z)\pi_{Z,\Theta}(z,\theta) \\
    &\propto& {\rm exp}\left(-\frac 12 \|\mC^{-1/2}(\mW z + \mu_{DCE} - y)\|^2 -\sum_{k=1}^L\frac{\|z^{(k)}\|^2}{2\theta_k} - \Phi(\theta)
    \right),
\end{eqnarray*}
where
\[
 \Phi(\theta) = \sum_{j=1}^L\left(\frac{\theta_j}{\vartheta_j}\right)^r + \sum_{j=1}^L\left(\frac{r_k}{2} +1 -r\beta_j \right)\log\theta_j.
\]
The structure of the posterior density is well suited  for the IAS algorithm, an efficient alternating scheme for computing a maximum a posteriori (MAP) estimate, extensively discussed in the literature \cite{calvetti2023ias}. The IAS algorithm can be summarized as follows:
\begin{enumerate}
    \item Initialize $\theta^0 = \vartheta$, setting the counter $\nu=0$.
    \item Repeat until a convergence criterion is met:
    \begin{enumerate}
        \item Given the current value $\theta^\nu$ update $z = z^{\nu+1}$ by minimizing
        \[
         E_\theta(z) = \frac 12 \|\mC^{-1/2}(\mW z + \mu_{DCE} - y)\|^2 +\sum_{k=1}^L\frac{\|z^{(k)}\|^2}{2\theta_k}, \quad \theta = \theta^\nu.
        \]
        \item Given the current value $z=z^{\nu+1}$ update $\theta = \theta^{\nu+1}$ by minimizing
        \[
         E_z(\theta) \sum_{k=1}^L\frac{\|z^{(k)}\|^2}{2\theta_k} + \Phi(\theta), \quad z = z^{\nu++1}
        \]
        \item Advance the counter $\nu\rightarrow \nu+1$ and check the convergence.
    \end{enumerate}
\end{enumerate}
The efficiency of the algorithm, compared, e.g., with $\ell^1$-based penalty functionals for promoting sparsity, arises from the fact that the updating of $z$ only requires the solution of a quadratic least squares problem, while the updating of $\theta$ can be performed component-wise, and the one-dimensional minimization problems for each component are simple to implement, in some cases admitting an explicit closed-form formula. We refer to \cite{calvetti2023bayesian} for general discussion of the IAS algorithm and of the parameters $\beta_k$ and $\vartheta_k$, and to  \cite{calvetti2025distributed,calvetti2025subspace} for recent computational developments. To overcome the non-convexity of the optimization problem when $r<1$, a hybrid version of the algorithm was proposed. The hybrid scheme starts the iteration with $r=1$ characterized by guaranteed convexity of the objective function and thus  convergence to a unique minimizer, and at some point switches to $r<1$ leading to fast local convergence with stronger sparsity promotion, see, e.g., \cite{calvetti2020sparsity}. As pointed out in the cited articles, the choices of the values of the other parameters $\beta_j$ and $\vartheta_j$ are related to the sparsity and sensitivity, and detailed out when the computed examples are described.  

The identification of the active brain regions corresponding to an observation can be based on monitoring the variance parameter $\theta$: The active regions are those associated with the most significant values, i.e.,  the regions for which
\begin{equation}\label{discard}
 \theta_j > p\,\max_{1\leq j\leq L}(\theta_j),
\end{equation}
where $p<1$ is a threshold parameter, e.g., $p =0.05$, are the most likely support of the activation underneath the MEG signal.

\subsection{Dictionary deflation}

The  Bayesian dictionary compression algorithm described above was originally proposed in \cite{bocchinfuso2024bayesian}. When its viability was tested on a hyperspectral data set, it was observed that the identification of the relevant subdictionary occasionally could be ambiguous: When the test data was arising from a single subdictionary, the correct class was typically not discarded by the variance thresholding criterion (\ref{discard}), however, the correct dictionary did not always correspond to the absolute maximum component of $\theta$, leading to a failure of the winner-takes-all classifier algorithm. This motivated the addition of a refining step to improve the identification success: After identifying a subset of relevant subdictionaries according to the criterion (\ref{discard}), the dictionary matching problem (\ref{dictionary matching}) is solved using a deflated dictionary $\mD^{\rm defl}$ obtained by removing form the original dictionary all subdictionaries deemed not relevant. The solution of the deflated dictionary matching problem is the minimizer of the objective function
\[
 E(x,\theta) = \|y - \mD^{\rm defl} x\|^2 +\sum_{j=1}^{n_d} \frac {x_j^2}{\theta_j} + \Psi(\theta),
\]
where $\Psi(\theta)$ accounts for the contribution from the Generalized gamma hyperprior model for the variance vector $\theta$. The minimizer can be computed effectively by the IAS algorithm with the promotion of sparsity for the individual components. As in Phase 1, the relevant dipoles are those with variance satisfying (\ref{discard}) and a winner-takes-all rule can be concentrate all activation to a single brain region. In the following section, we test systematically the performance of the proposed two-phase approach on MEG data with the subdictionaries defined by the Destrieux brain atlas.

\section{Assessment by simulations}

In this section we test how well the brain regions defined by the Destrieux atlas can be identified from MEG data with the proposed Bayesian dictionary learning methodology by means of the following numerical protocol. In addition to assessing the success of the proposed algorithms as a classifier, the results of the simulations open a discussion on whether  brain atlases based on anatomical brain features should be replaced by brain atlas focusing on the cerebral electromagnetic activity for studies based on MEG and EEG data, e.g., by identifying brain regions that are unambiguously identifiable by the available data.

We download the head model and the corresponding lead field matrices  $\mL_k \in\R^{m\times 3}$, $1\leq k\leq n$ from the Open MEG Archive (OMEGA) website {\tt
https://www.mcgill.ca/bic/neuroinformatics/omega},  see \cite{niso2016omega}. The discretized head model consist of $n = 15\, 002$ dipoles, and $m=306$ measurements performed by 102 sensors, each comprising a magnetometer and a gradiometer, recording the magnetic flux density as well as tangential gradient components. 
The source space is divided into $L = 148$ parcels according to the Destrieux brain atlas \cite{destrieux2010automatic} using the automatic parcellation provided by Brainstorm \cite{tadel2011brainstorm}, constituting the basis for the partition into subdictionaries.

\begin{figure}[h!]
\centerline{
\includegraphics[width=10cm]{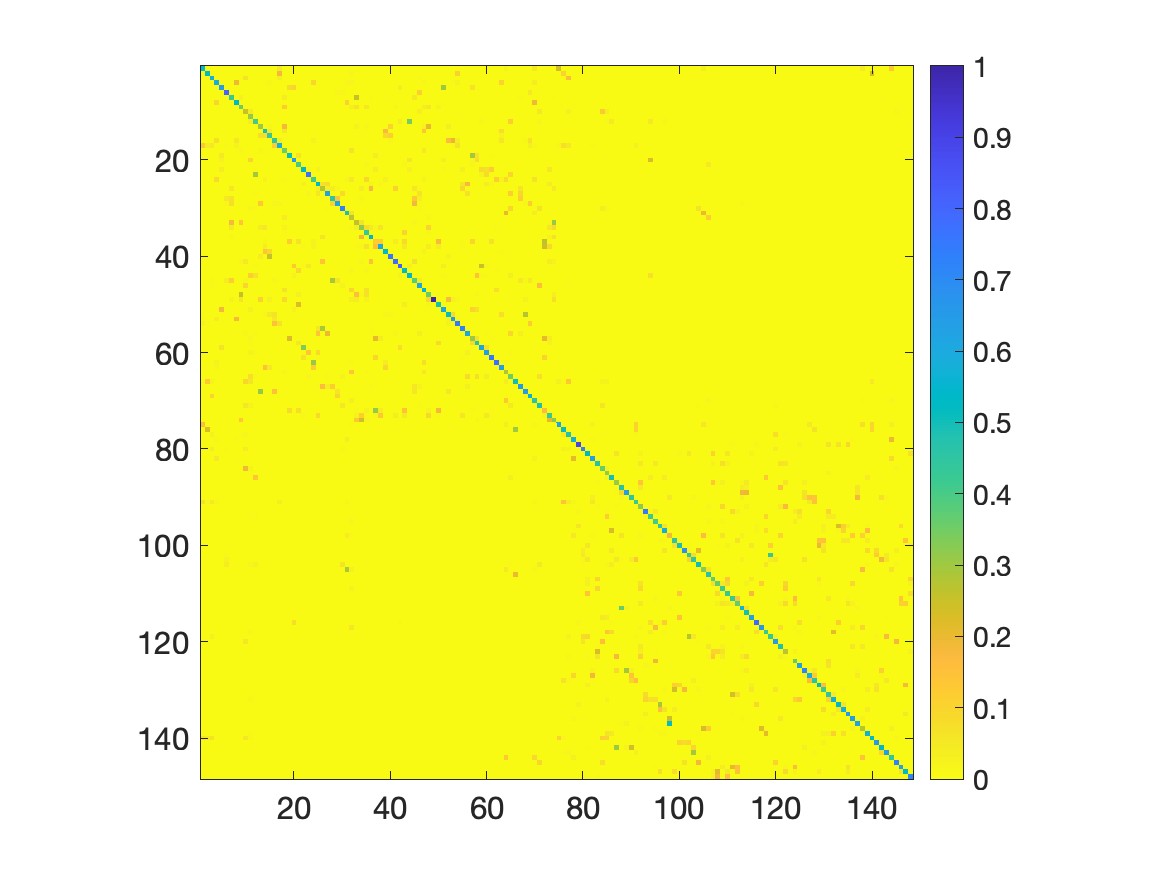}
\includegraphics[width=10cm]{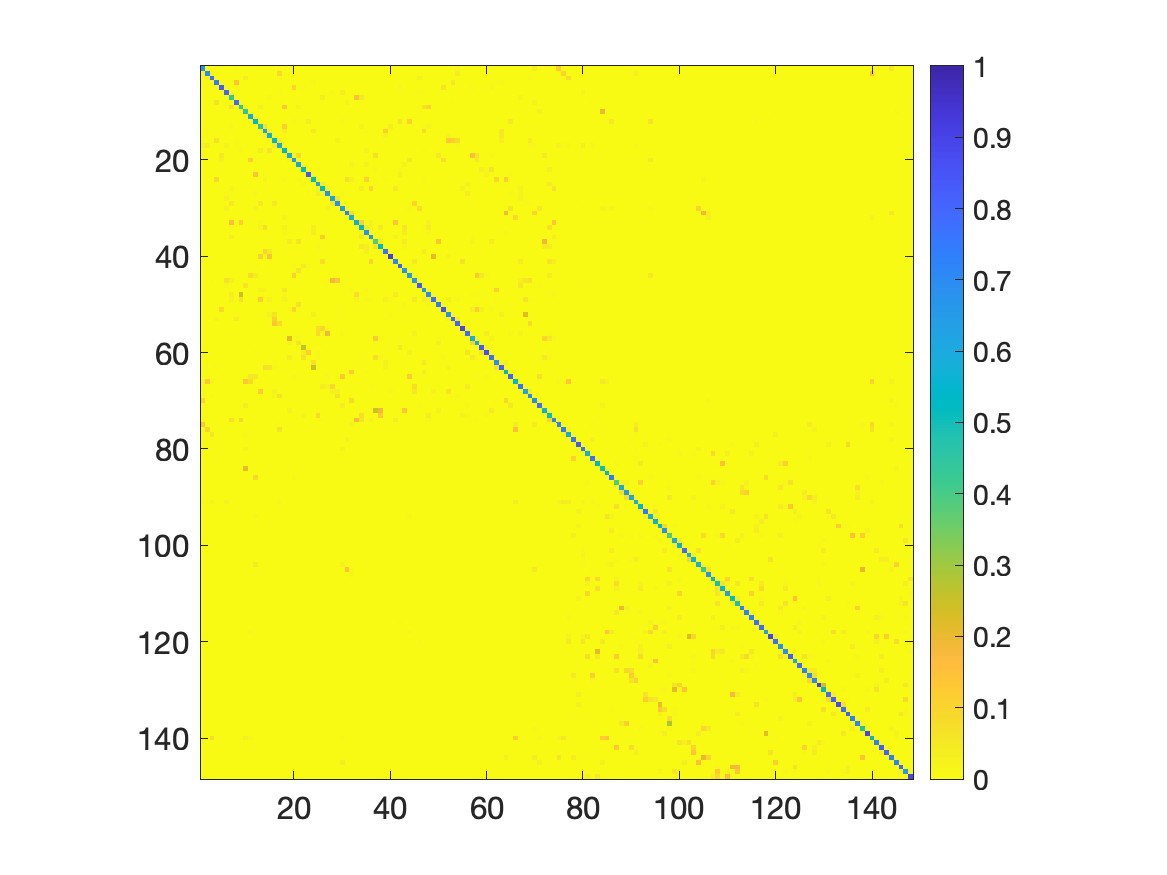}
}
\caption{Heatmaps of the confusion matrices after the first phase using reduced dictionaries and winner-takes-all classification principle (left) and second phase based on the deflated model (right).}\label{fig:confusion matrices}    
\end{figure}

For each brain region, we simulate an active patch of six adjacent dipoles at a randomly selected location, with  dipole amplitudes assigned randomly from the uniform [0,1] distribution, and compute the induced magnetic field at the magnetometers. The noiseless computed data are then corrupted by a scaled Gaussian white noise with standard deviation 0.5\% of the maximum noiseless component of the data vector. Starting with the simulated noisy MEG data, we  run the two-phase classifier as described in \cite{bocchinfuso2024bayesian}: The first phase is based on the hybrid group sparsity IAS algorithm, in which the MAP minimization for the compressed dictionary is first estimated using the gamma hyperprior ($r=1$), then switch to generalized gamma hyperprior model ($r=1/2$) to more strongly promote sparsity in the solution. At the end of Phase 1, we discard all dipole of the full dictionary in regions not satisfying (\ref{discard}) value $p = 0.005$, and Phase 2 solves the deflated dictionary matching problem using the IAS algorithm with the gamma hyperprior model ($r=1$).  We repeat this process 100 times for each brain region, for a total of 14\,800 classification tests, tallying the winner-takes-all classification results separately for Phase 1 and Phase 2 and storing the results in the confusion matrices $\mC^1$ and $\mC^2$, respectively.
We compute the column-scaled matrices $\mP^1$ and $\mP^2$, and from their heatmaps, shown in Figure~\ref{fig:confusion matrices}. It can be seen that, overall, the misclassification error is smaller at the end of Phase 2 than at the end of Phase 1. Additional information about the misclassification is provided by the other impurity measures introduced in section~\ref{sec:assessing}.

\begin{figure}
\centerline{
\includegraphics[width=10cm]{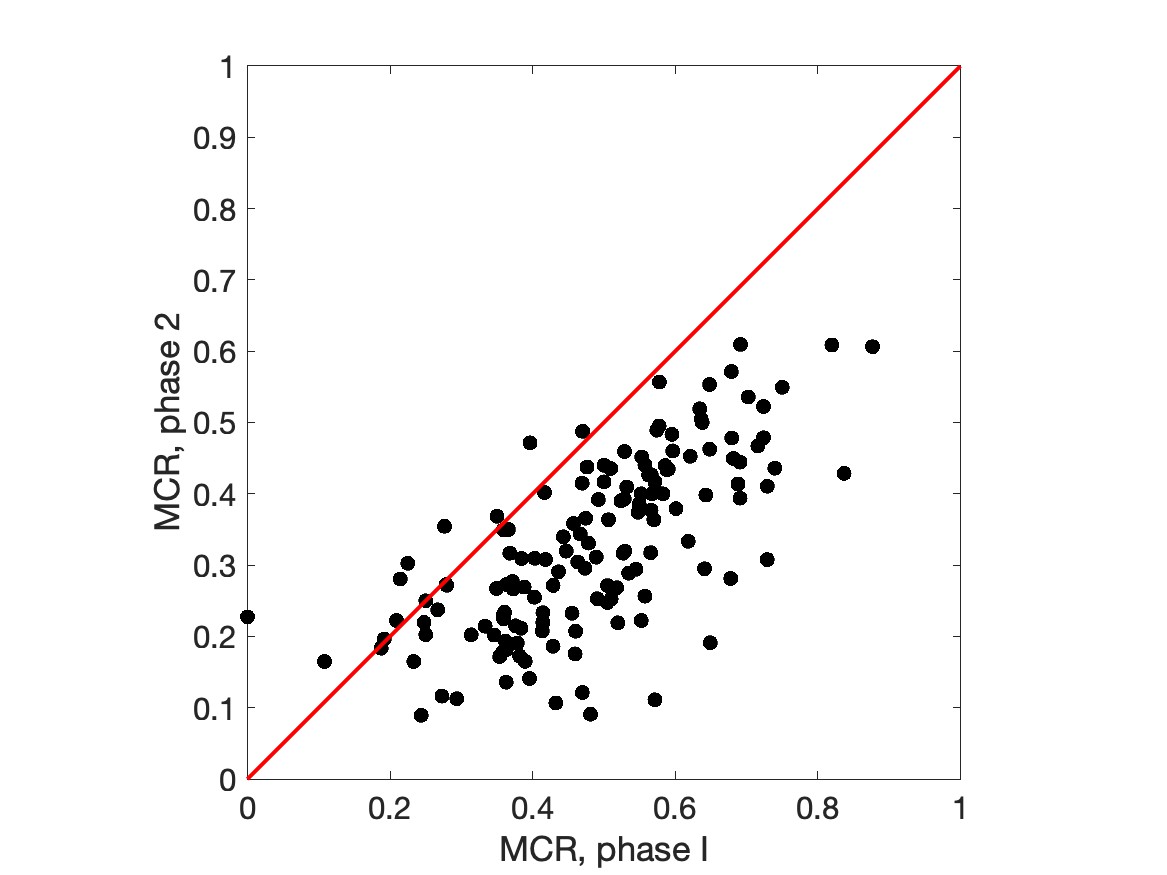}
\includegraphics[width=10cm]{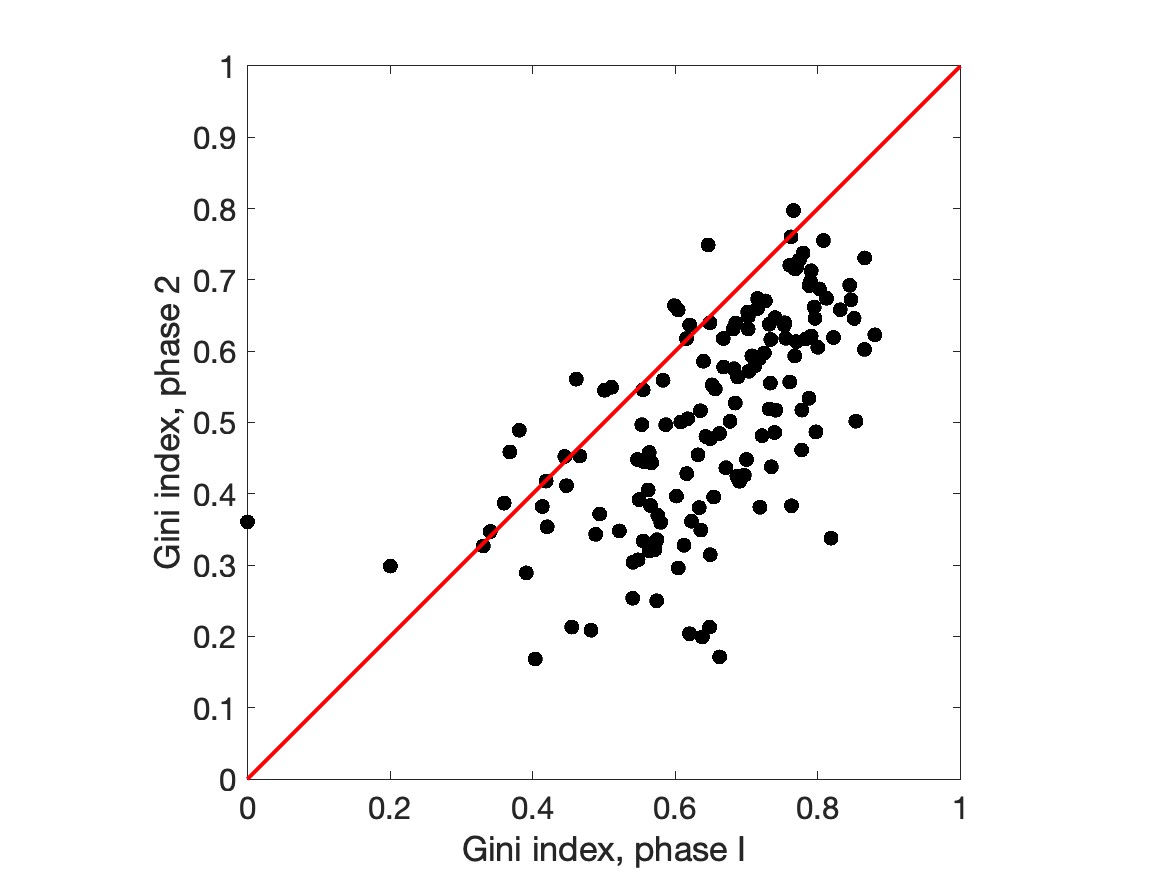}
}
\caption{Scatter plots of the misclassification errors of each brain region (left) and the Gini indices (right). The horizontal axes correspond to the classification based on the reduced dictionaries, and the vertical plot to the deflated dictionaries. Observe that decreases of the impurity in Phase 2 is indicated by the marker for brain region to be below the red diagonal line.}\label{fig:impurities}    
\end{figure}

Figure~\ref{fig:impurities} shows the misclassification rates $MCR_j$ and Gini indices $G_j$ for each brain region in scatter plot form, the horizontal axis corresponding to the result at the end of Phase I, and the vertical axis corresponding to the result at the end of Phase II. Dots below the main diagonal indicate that for the corresponding brain regions the classification improves in Phase II, as we would expect and as is indeed the case for most of the brain regions, although in some rare case the impurity is higher at the end of Phase 2 than at the end of Phase 1.

\begin{figure}
\centerline{\includegraphics[width=1.1\textwidth]{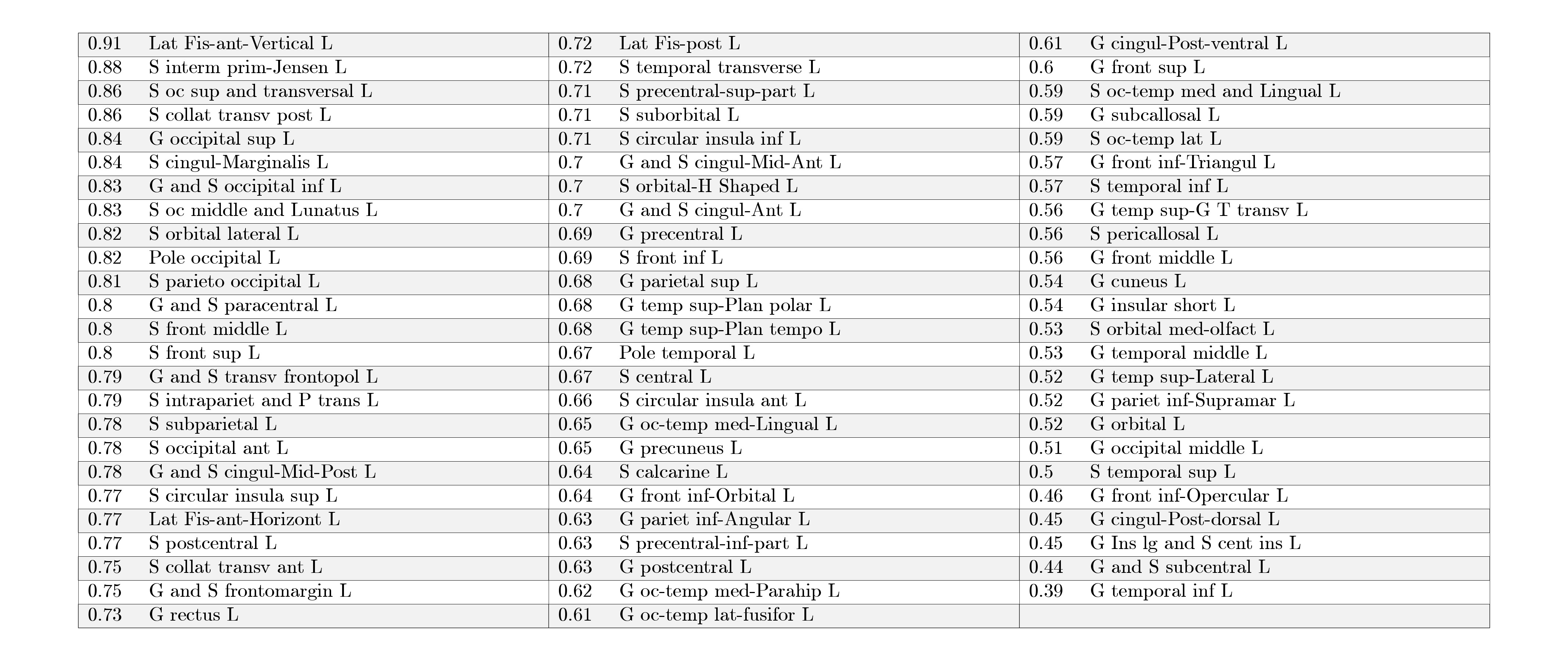}}
\caption{The recalls $s_j = 1-MCR_j$ of the brain regions of the left hemisphere in decreasing order. Here, both phases of the classifier algorithm was applied. }\label{fig:recalls}
\end{figure}

Following the terminology used in binary classification, we define the recalls of each brain region as
\[
 s_j = 1 - MCR_j.
\]
In Figure~\ref{fig:recalls}, the 74 Destrieux regions in the left hemisphere are listed in decreasing  order of recall computed at the end of Phase 2. 

To gain some understanding on why some brain regions may be prone to misclassified and to shed some light on which brain regions the confusion entails, we investigate in more depth four regions with different classification success rate, and visualize the corresponding column vectors of the matrix $\mP$, recalling that the $i$th entry of the column can be viewed as the probability of the activity in region $i$ to be classified as activity in the region defining the column. If the classifier is reliable, we would expect the support of the column to be concentrated in the $i$th entry.
. 

\begin{figure}[h]
\centerline{
\includegraphics[width=\textwidth]{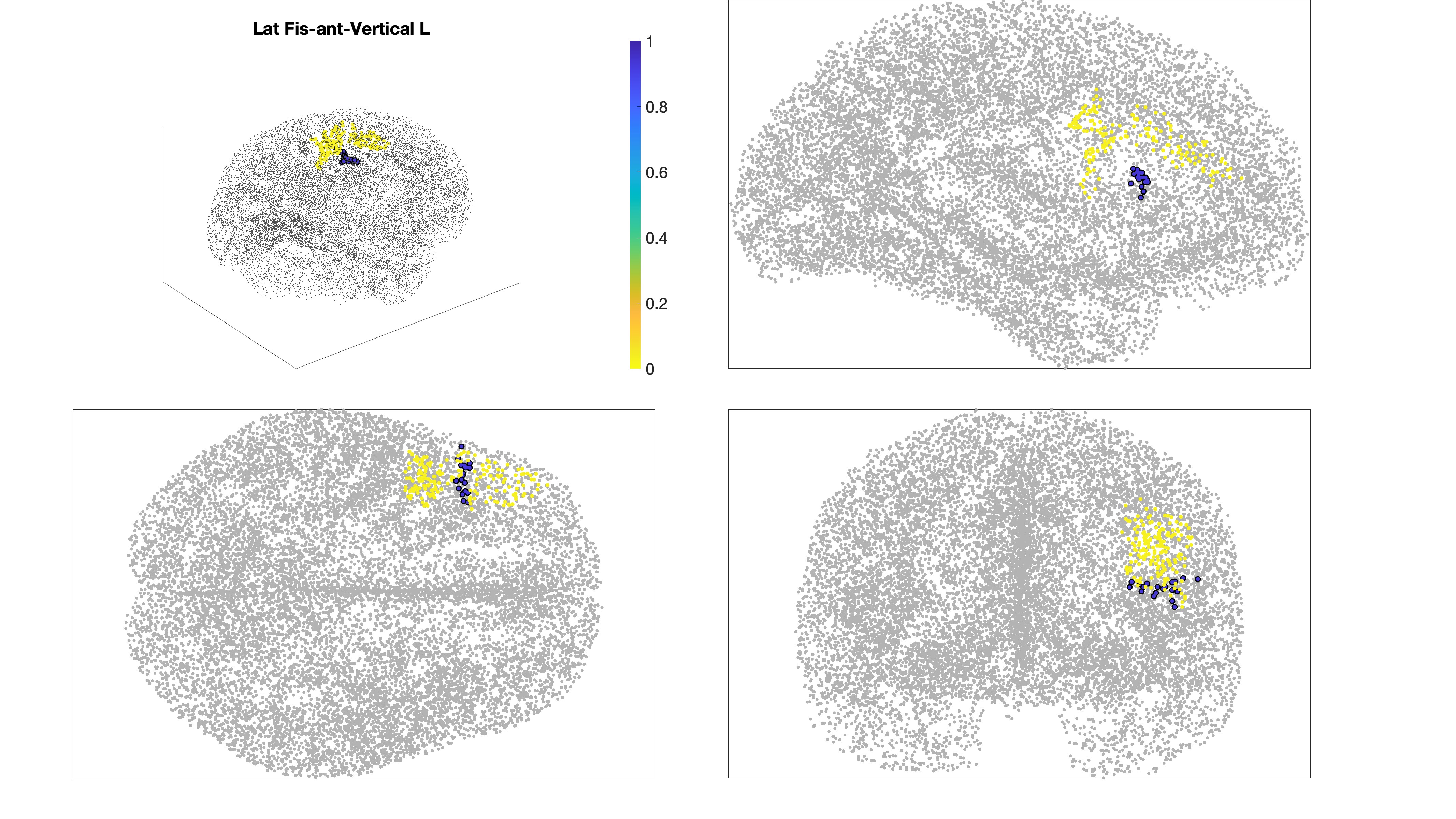}
}
\caption{Visualization of entries  $(\mP_{ij})_{i=1}^N$, $j$ of the $j$th vector, interpreted as probabilities, corresponding to vertical ramus of the anterior segment of the lateral fissure, the brain region with the highest recall value of 0.91.}\label{fig:lat_fiss}
\end{figure}

\begin{figure}
\centerline{
\includegraphics[width=\textwidth]{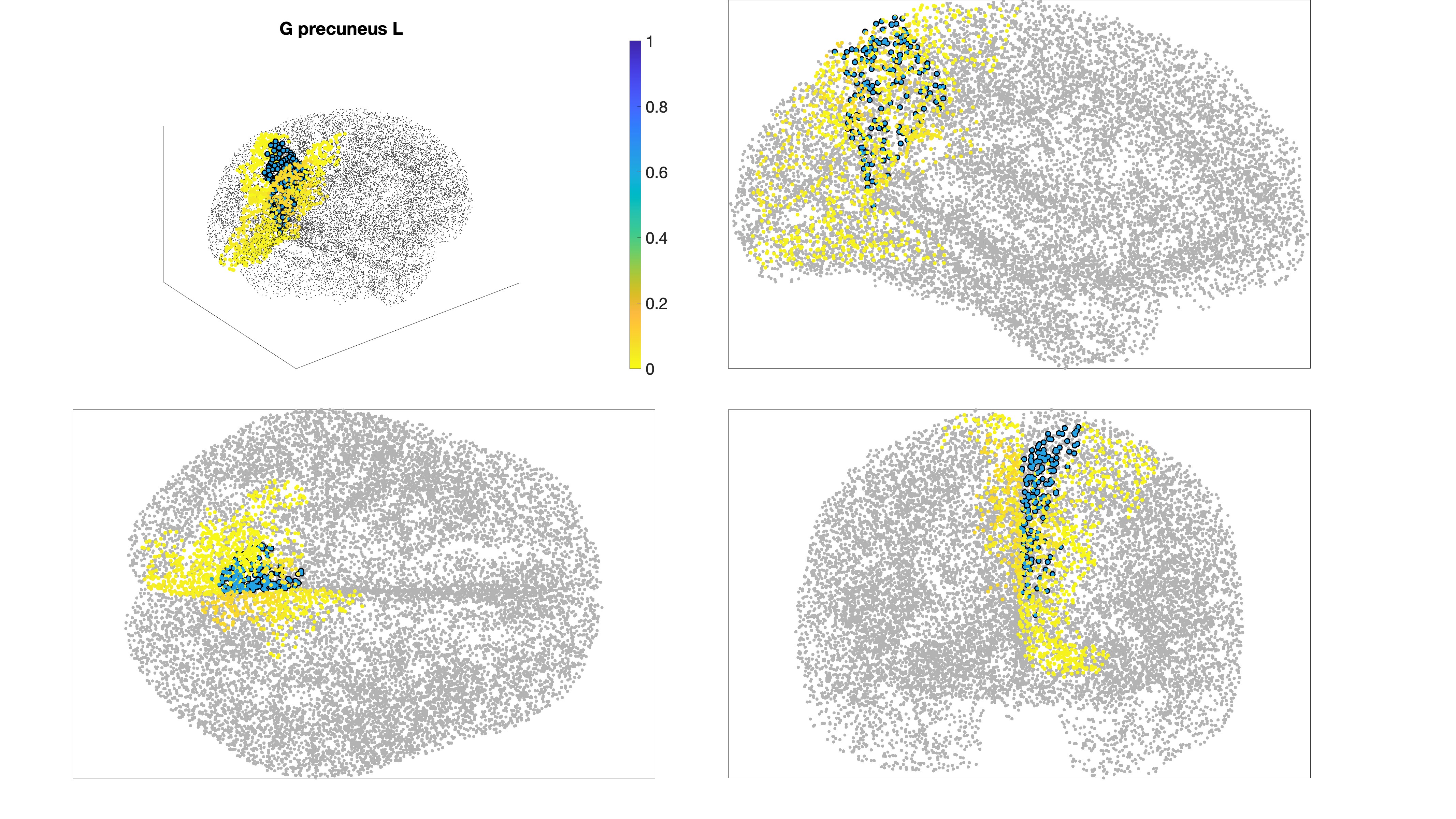}
}
\caption{Visualization of the probabilities of the regions to be identified as left precuneus, with recall value 0.65. Observe that due to the proximity of the left precuneus to the longitudinal fissure separating the hemispheres, the misclassification includes some brain regions on the right hemisphere.}\label{fig:precuneus}
\end{figure}

\begin{figure}
\centerline{
\includegraphics[width=\textwidth]{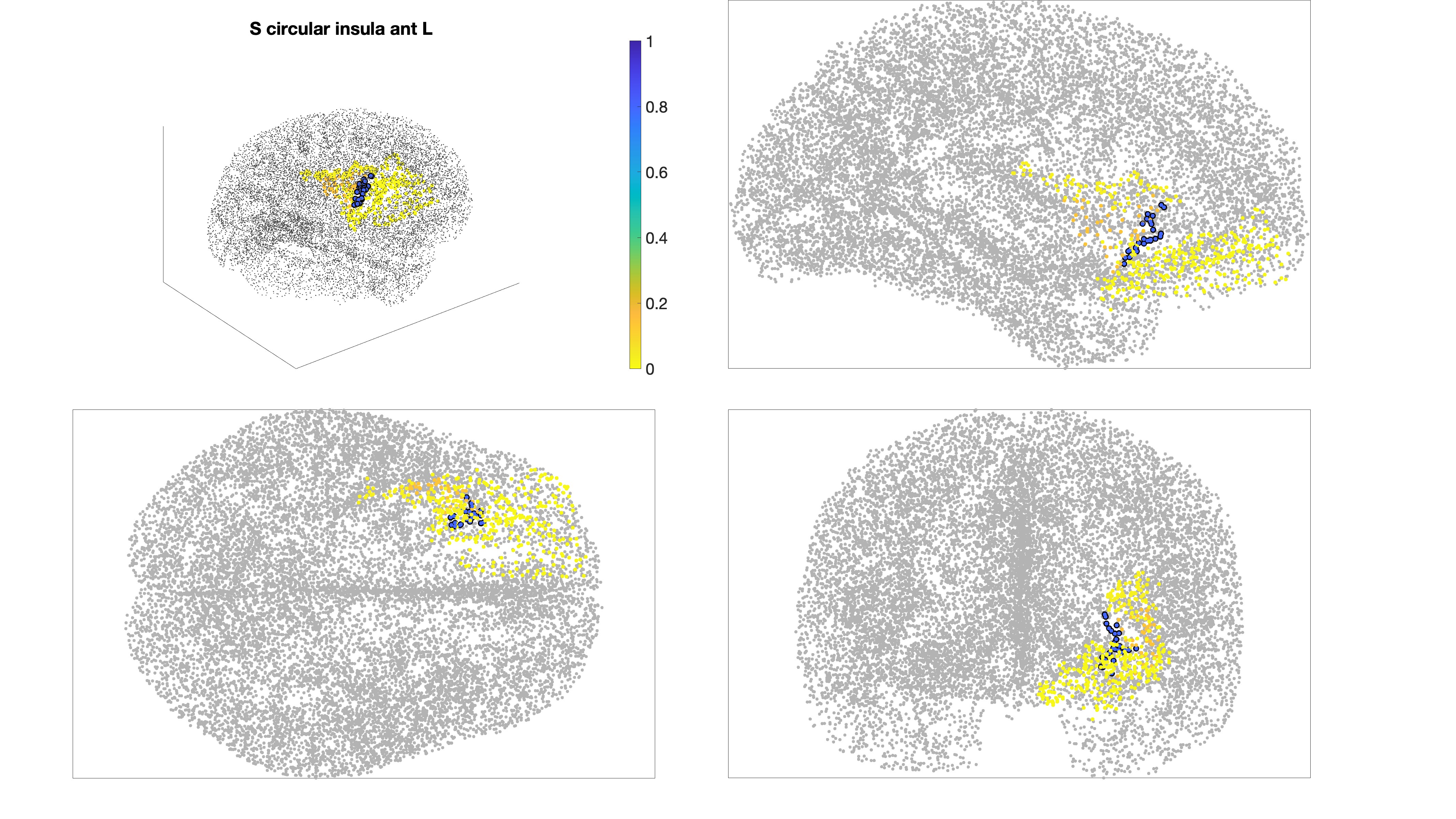}
}
\caption{Visualization of the probabilities of brain regions to be identified as the anterior sulcus of the left circular insula. This deep brain region has the recall value 0.66.}\label{fig:insula}
\end{figure}

\begin{figure}
\centerline{
\includegraphics[width=\textwidth]{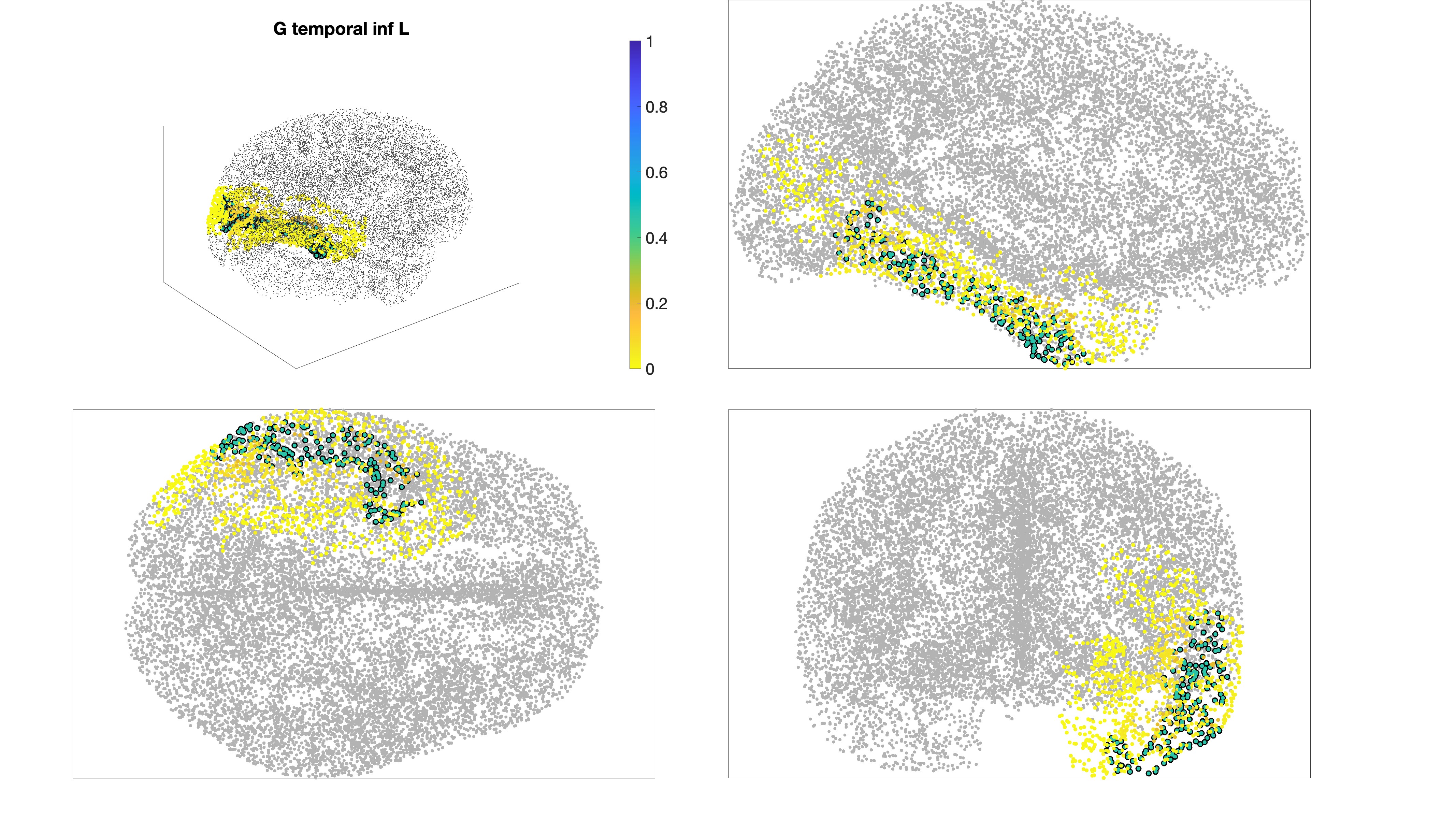}
}
\caption{Visualization of the probabilities of brain regions to be identified as the left inferior temporal gyrus, the region with the lowest recall, 0.39, in the left hemisphere. The confusion may be due to the relatively large size of the region, as well as its position in the brain making the MEG data less sensitive to signals arising from it.}\label{fig:temporal}
\end{figure}

Figure~\ref{fig:lat_fiss} refers to the anterior segment of the lateral fissure in the left hemisphere, the brain region having highest recall. Figures \ref{fig:precuneus} and \ref{fig:insula} correspond to the left precuneus and to the anterior sulcus of the left circular insula, two brain regions with mid-range recall values, and Figure \ref{fig:temporal} to the inferior temporal gyrus, the brain region in the left hemisphere with the lowest recall. One reason for the low recall for the inferior temporal gyrus is its distance from the MEG sensors, contributing to low sensitivity of the MEG data to signals coming from this region. Moreover, the relatively large extension of this brain region contributes to a wide variability in the induced magnetic field, increasing the probability that the signal is mistakenly identified as coming from one of the numerous neighboring regions. The regions with the highest recalls, on the other hand, tend to be more superficial, close to the MEG the sensors, and relatively small. While most of the confusion is confined to the hemisphere where the identified region lies, we see in Figure~\ref{fig:precuneus}
that the proximity to the longitudinal fissure may spread the confusion across the hemispheres.

\subsection{The MEG identification trees}

The column-normalized confusion matrix $\mP$ suggests an alternative interpretation of the probability of the activity attributed to a given ROI as a directed network. More specifically, each column of $\mP$ can be associated with a weighted directed tree graph with a self-loop at the root, which is the node associated with the ROI indexing the column. The remaining nodes are the other ROIs in the atlas and the entries of the column are the weights of the edges from the root node. The complexity of the trees depends on how well the root ROI is identifiable: for most ROIs the tree has only a few nodes, and only very few ROIs are associated with a tree with many leaves.

Recalling that anatomical considerations have been the guiding principle behind the Destrieux atlas, and in view of the tight packing of gyri and sulci, it is not unreasonable that some of ROIs are very difficult to identify from MEG data. The location of a ROI relative to the placement of the sensors measuring the magnetic field, its volume and its positions in relation  to its neighbors are factors playing an important role in how well it can be recovered. 

We illustrate the tree interpretation for four ROIs in the left hemisphere from the Destrieux atlas, the \emph{middle and posterior cingulate region (sulcus and gyrus)}, the \emph{frontal superior gyrus}, the \emph{inferior parietal angular gyrus} and the \emph{precuneus gyrus}. The selected brain regions play a central role in the default mode network (DMN), a focus of interest in the brain connectivity studies \cite{raichle2015brain,menon202320}.

\begin{figure}
\centerline{
\includegraphics[width=\textwidth]{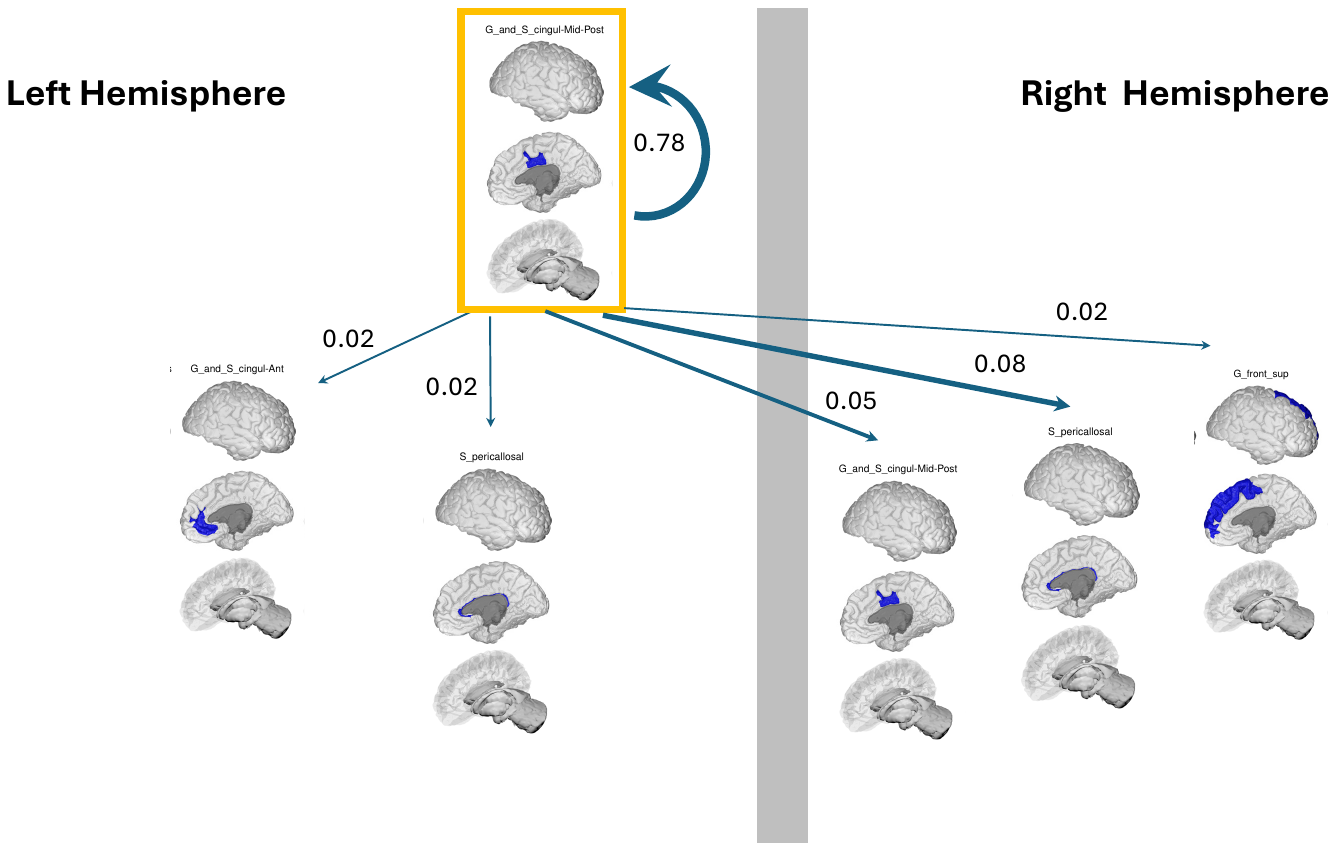}
}
\caption{\label{fig:MidPosteriorGraph} Graph representation of the attribution of activity restricted to the left mid-posterior cingulate region from MEG data.}\end{figure}

The left middle and posterior cingulate cortex and sulcus (MPCCS) play a role in some important higher-level functions, including allocation of attention, anticipation of rewards, decision-making, control of impulses, and emotions. The weighted tree associated with the left MPCCS region, shown in figure~\ref{fig:MidPosteriorGraph}, suggests that cerebral activity localized to this brain region is very likely to be correctly and unequivocally identified from MEG measurement, as it was the case in 78\% of the test cases. The anatomical location of this region, right above the pericallosal sulcus may suggest a challenge in the correct distinction between the two areas, however in our simulations only 10\% of the times activity from the left MPCCS region was localized in the left or right pericallosal sulcus. The closeness of the MPCCS to the to midsagittal plane explains why 5\% of the cases, the activity was attributed to the specular, and very near right mid-posterior cingulate region, which in our simulations occurred in only 5\% of the cases. 

\begin{figure}
\centerline{
\includegraphics[width=\textwidth]{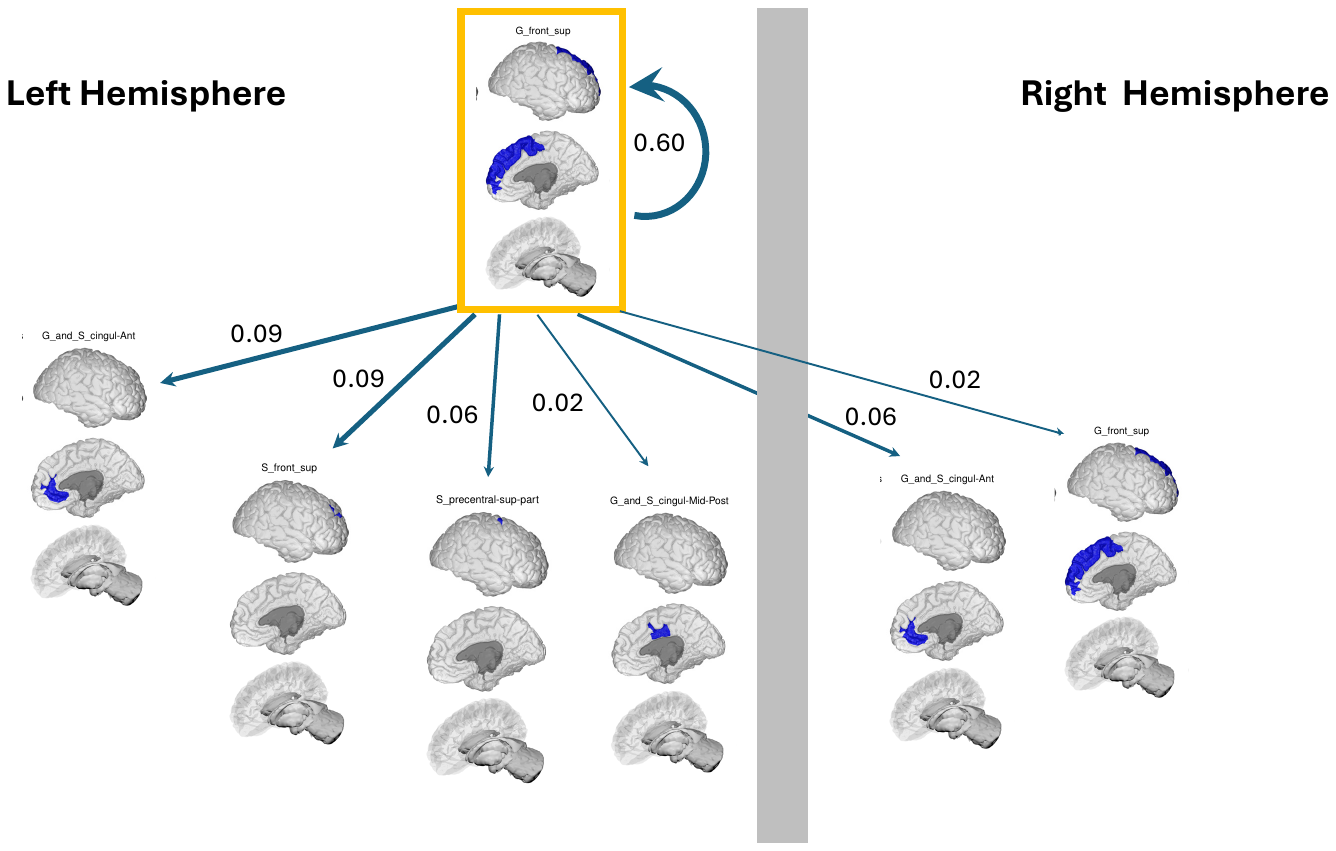}
}
\caption{\label{fig:FrontalSuperiorGraph} Graph representation of the attribution of activity originating inside the left frontal superior gyrus  from MEG data.}
\end{figure}

The  left superior frontal gyrus (SFG), which is the top part of the prefrontal cortex, is a node of the cognitive control network and the default mode network. It has been hypothesized that the SFG, which is primarily associated with Brodmann areas 6, 8, and 9, plays a role in higher cognitive functions, for example working memory. 
self-monitoring, planning, and organization. 
The weighted tree associated with the left SFG  region, shown in figure~\ref{fig:FrontalSuperiorGraph}, indicates that the most likely regions to be confused with this ROI are, not surprisingly, the  left superior frontal sulcus, the left and right cingulate anterior sulcus and gyrus, and the left precentral superior parietal sulcus, which are all anatomically very close. 

The left angular gyrus (LAG) is a horseshoe shaped region located in the inferior portion of the parietal lobule, involved with involved in various cognitive functions, including language processing, spatial cognition, and memory retrieval. The identification of activity in this region from MEG signal may be attributed to the superior temporal, intraparietal and Jensen sulci, which are anatomically very close and very small, and, with a lower probability, to other regions within the parietal lobe. 

The left precuneus gyrus (LPG) is located in the portion of the superior parietal lobule on the medial surface, in front of the upper portion of the occipital lobe. LPG is bounded by the parieto-occipital sulcus and by the subparietal sulcus, thus it is not surprising that these are the more frequent confounding regions in source reconstruction from MEG.   

To summarize the findings, we note that the Destrieux parcellation with its distinction between sulci and gyri, based on widely accepted anatomical conventions, has been designed for the automatically labeling the whole cerebral cortex, but may be suboptimal for assessing the resolution power of MEG data for source identification. 

\begin{figure}
\centerline{
\includegraphics[width=\textwidth]{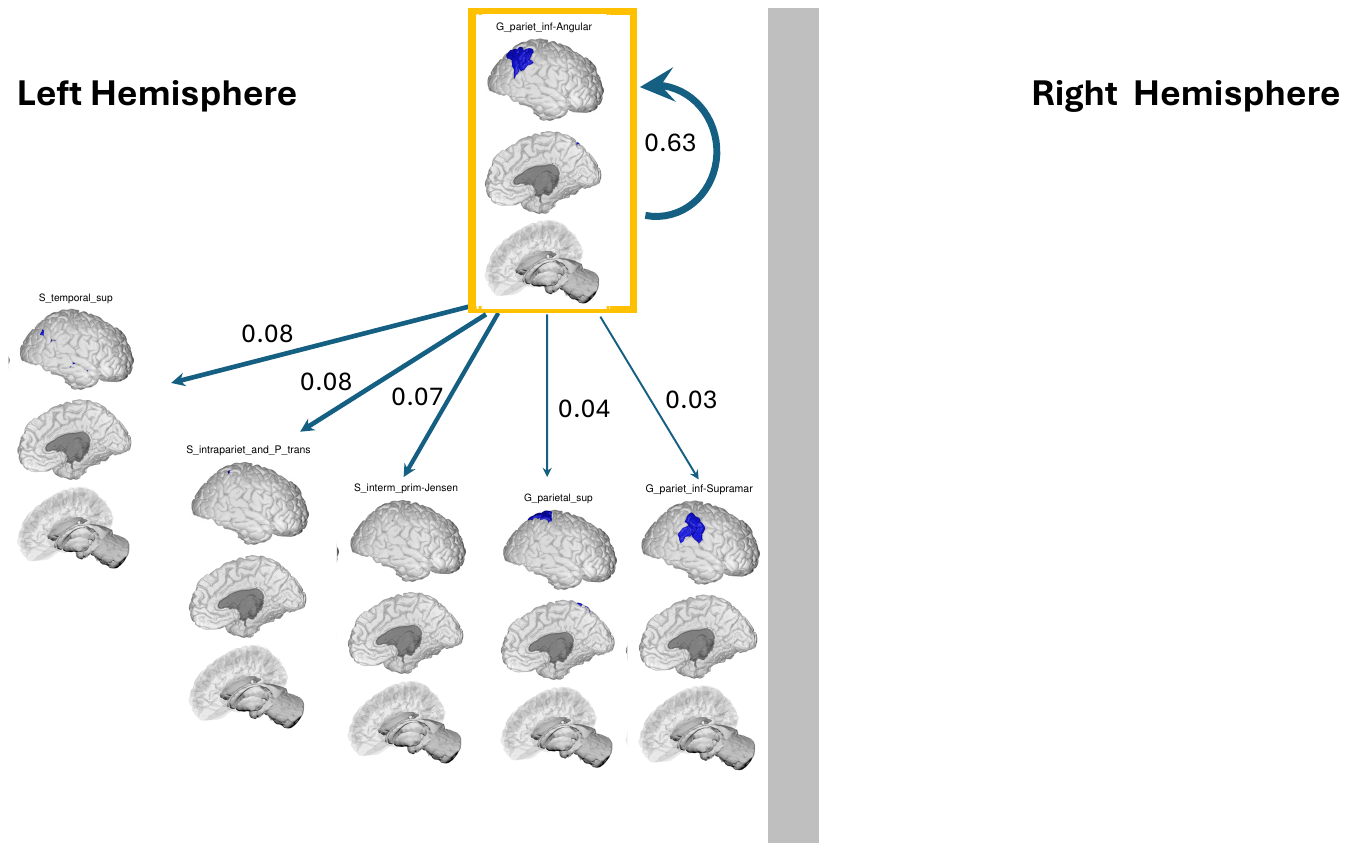}
}
\caption{\label{fig:ParietallAngularGraph}Graph representation of the attribution of activity originating inside the left angular gyrus in the parietal lobe from MEG data. }
\end{figure}

\begin{figure}
\centerline{
\includegraphics[width=\textwidth]{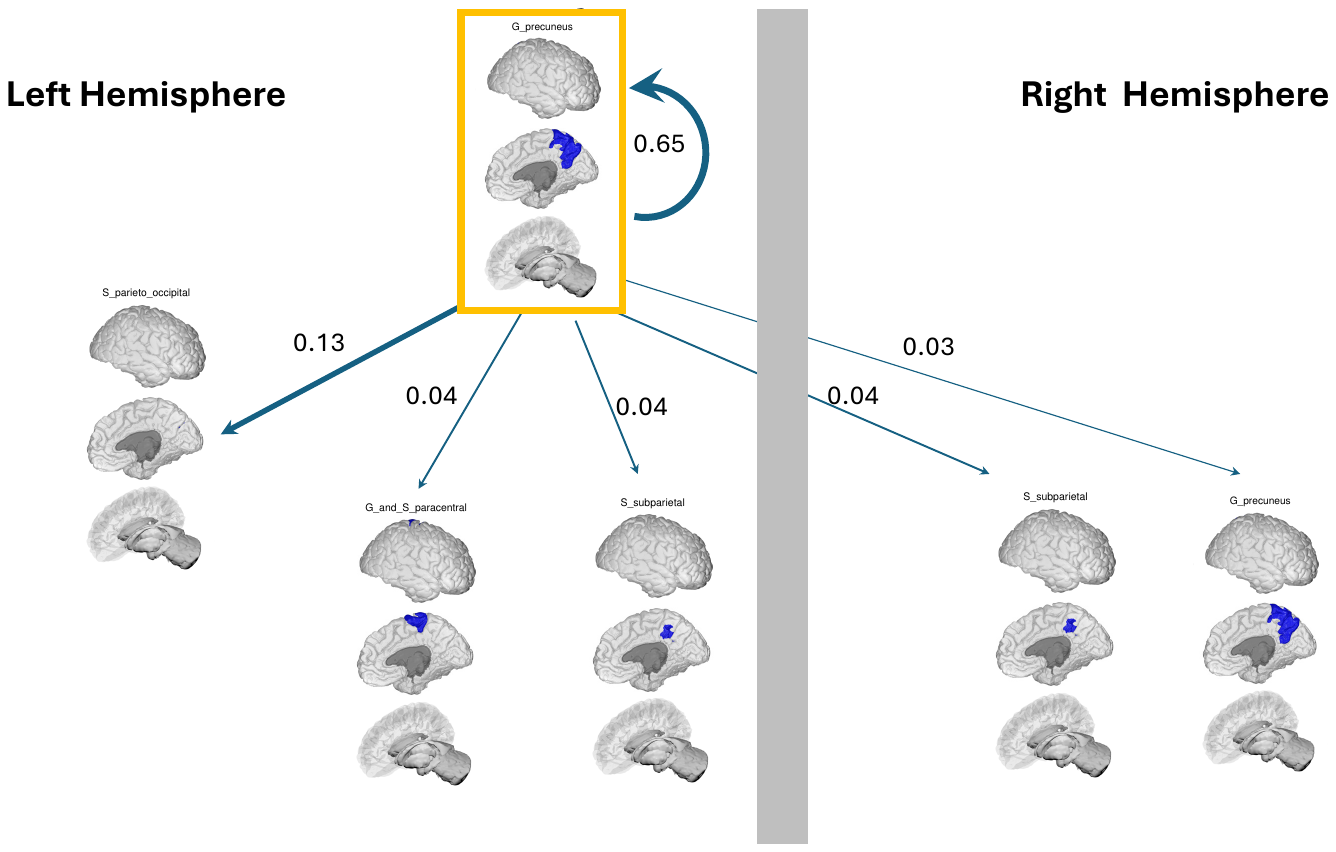}
}
\caption{\label{fig:PrecuneusGraph}Graph representation of the attribution of activity originating inside the left precuneus gyrus i from MEG data.}
\end{figure}

\section{Discussion}

In this article, a version of the Bayesian dictionary learning algorithm was adapted for the MEG source identification problem. The identification of each region  of the the Destrieux atlas from MEG data was systematically tested through numerical simulations by generating a random patch activity in a single brain region at the time. The performance of the classifier was assessed over a suite of simulated data arising from 100 randomly activated patches for each brain region.  From the point of view of the practical applications, the question was posed as an a posteriori estimate: What are the probabilities of various regions to be active, given that the winner-takes-all algorithm identifies one region? While in this question, it is implicitly assumed that the activity is confined in one single region, the analysis sheds light on the success of identifying more than one activity region from the data: The the testing method was posed only for having an easily quantifiable output, but it is not an inherent feature of the algorithm itself, as the algorithm does not require the winner-takes-all identification. Indeed, the output of the algorithm can be summarized in terms of the estimated variance vectors $\theta$, and by thresholding, several brain regions can be interpreted to be simultaneously active.

The Bayesian formulation underlying the dictionary learning algorithm opens possibilities for uncertainty quantification beyond the dictionary identification of the active region. As pointed out in \cite{bocchinfuso2024bayesian}, the first phase of the algorithm selects a few brain regions and while usually the activated one is among them, it is not necessarily the one with the highest variance value, thus the winner-take-all criterion would lead to misclassification, which is the main motivation for the addition of Phase II.  Since it may be argued that the solution at the end of Phase I is a maximum a posteriori (MAP) estimate, which may be seen as not a true Bayesian solution because of its dependence on the model parametrization, an alternative approach is to find instead the  computationally more expensive sampling-based posterior mean estimate. In \cite{calvetti2025subspace}, this question was addressed by means of developing a fast MCMC sampler suitable for problems with priors favoring  group sparsity, and it was shown that the current model for MEG is amenable for this approach. While the sampling based method may not be as suitable for processing large time series data as the current dictionary learning method, it provides a way to assess the reliability of the algorithm at selected control points. 

\section*{Acknowledgements}

The work of DC was partly supported by the NSF grants DMS 1951446 and  DMS 2513481, and that of  ES by the NSF grant DMS 2204618 and DMS 2513481.

\bibliographystyle{siam} 
\bibliography{biblio}

\begin{thebibliography}{10}

\bibitem{aharon2006k}
{\sc M.~Aharon, M.~Elad, and A.~Bruckstein}, {\em K-svd: An algorithm for designing overcomplete dictionaries for sparse representation}, IEEE Transactions on signal processing, 54 (2006), pp.~4311--4322.

\bibitem{baillet2001electromagnetic}
{\sc S.~Baillet, J.~C. Mosher, and R.~M. Leahy}, {\em Electromagnetic brain mapping}, IEEE Signal processing magazine, 18 (2001), pp.~14--30.

\bibitem{bocchinfuso2024bayesian}
{\sc A.~Bocchinfuso, D.~Calvetti, and E.~Somersalo}, {\em Bayesian sparsity and class sparsity priors for dictionary learning and coding}, Journal of Computational Mathematics and Data Science, 11 (2024), p.~100094.

\bibitem{calvetti2015hierarchical}
{\sc D.~Calvetti, A.~Pascarella, F.~Pitolli, E.~Somersalo, and B.~Vantaggi}, {\em {A hierarchical {K}rylov--{B}ayes iterative inverse solver for {MEG} with physiological preconditioning}}, Inverse Problems, 31 (2015), p.~125005.

\bibitem{calvetti2023ias}
{\sc D.~Calvetti, A.~Pascarella, F.~Pitolli, E.~Somersalo, and B.~Vantaggi}, {\em The {IAS}-{MEEG} package: a flexible inverse source reconstruction platform for reconstruction and visualization of brain activity from m/eeg data}, Brain Topography, 36 (2023), pp.~10--22.

\bibitem{calvetti2020sparsity}
{\sc D.~Calvetti, M.~Pragliola, and E.~Somersalo}, {\em Sparsity promoting hybrid solvers for hierarchical bayesian inverse problems}, SIAM Journal on Scientific Computing, 42 (2020), pp.~A3761--A3784.

\bibitem{calvetti2023bayesian}
{\sc D.~Calvetti and E.~Somersalo}, {\em Bayesian Scientific Computing}, vol.~215, Springer Nature, 2023.

\bibitem{calvetti2025distributed}
\leavevmode\vrule height 2pt depth -1.6pt width 23pt, {\em Distributed tikhonov regularization for ill-posed inverse problems from a bayesian perspective}, Computational Optimization and Applications,  (2025), pp.~1--32.

\bibitem{calvetti2025subspace}
\leavevmode\vrule height 2pt depth -1.6pt width 23pt, {\em Subspace splitting fast sampling from gaussian posterior distributions of linear inverse problems}, arXiv preprint arXiv:2502.05703,  (2025).

\bibitem{calvetti2019hierachical}
{\sc D.~Calvetti, E.~Somersalo, and A.~Strang}, {\em Hierachical {B}ayesian models and sparsity: $\ell^2$-magic}, Inverse Problems, 35 (2019), p.~035003.

\bibitem{desikan2006automated}
{\sc R.~S. Desikan, F.~S{\'e}gonne, B.~Fischl, B.~T. Quinn, B.~C. Dickerson, D.~Blacker, R.~L. Buckner, A.~M. Dale, R.~P. Maguire, B.~T. Hyman, et~al.}, {\em An automated labeling system for subdividing the human cerebral cortex on mri scans into gyral based regions of interest}, Neuroimage, 31 (2006), pp.~968--980.

\bibitem{destrieux2010automatic}
{\sc C.~Destrieux, B.~Fischl, A.~Dale, and E.~Halgren}, {\em Automatic parcellation of human cortical gyri and sulci using standard anatomical nomenclature}, Neuroimage, 53 (2010), pp.~1--15.

\bibitem{doucet2011brain}
{\sc G.~Doucet, M.~Naveau, L.~Petit, N.~Delcroix, L.~Zago, F.~Crivello, G.~Jobard, N.~Tzourio-Mazoyer, B.~Mazoyer, E.~Mellet, et~al.}, {\em Brain activity at rest: a multiscale hierarchical functional organization}, Journal of neurophysiology, 105 (2011), pp.~2753--2763.

\bibitem{gillis2020nonnegative}
{\sc N.~Gillis}, {\em Nonnegative matrix factorization}, SIAM, 2020.

\bibitem{hamalainen1993magnetoencephalography}
{\sc M.~H{\"a}m{\"a}l{\"a}inen, R.~Hari, R.~J. Ilmoniemi, J.~Knuutila, and O.~V. Lounasmaa}, {\em Magnetoencephalography—theory, instrumentation, and applications to noninvasive studies of the working human brain}, Reviews of modern Physics, 65 (1993), p.~413.

\bibitem{hilgetag2020hierarchy}
{\sc C.~C. Hilgetag and A.~Goulas}, {\em ‘{Hierarchy}’ in the organization of brain networks}, Philosophical Transactions of the Royal Society B, 375 (2020), p.~20190319.

\bibitem{ilmoniemi2019brain}
{\sc R.~J. Ilmoniemi and J.~Sarvas}, {\em Brain signals: physics and mathematics of {MEG} and {EEG}}, MIT Press, 2019.

\bibitem{kaipio2006statistical}
{\sc J.~Kaipio and E.~Somersalo}, {\em Statistical and computational inverse problems}, vol.~160, Springer Science \& Business Media, 2006.

\bibitem{kaipio2007statistical}
\leavevmode\vrule height 2pt depth -1.6pt width 23pt, {\em Statistical inverse problems: discretization, model reduction and inverse crimes}, Journal of computational and applied mathematics, 198 (2007), pp.~493--504.

\bibitem{krishnaswamy2017sparsity}
{\sc P.~Krishnaswamy, G.~Obregon-Henao, J.~Ahveninen, S.~Khan, B.~Babadi, J.~E. Iglesias, M.~S. H{\"a}m{\"a}l{\"a}inen, and P.~L. Purdon}, {\em Sparsity enables estimation of both subcortical and cortical activity from meg and eeg}, Proceedings of the National Academy of Sciences, 114 (2017), pp.~E10465--E10474.

\bibitem{lancaster2000automated}
{\sc J.~L. Lancaster, M.~G. Woldorff, L.~M. Parsons, M.~Liotti, C.~S. Freitas, L.~Rainey, P.~V. Kochunov, D.~Nickerson, S.~A. Mikiten, and P.~T. Fox}, {\em Automated talairach atlas labels for functional brain mapping}, Human brain mapping, 10 (2000), pp.~120--131.

\bibitem{lim2017sparse}
{\sc M.~Lim, J.~M. Ales, B.~R. Cottereau, T.~Hastie, and A.~M. Norcia}, {\em Sparse eeg/meg source estimation via a group lasso}, PloS one, 12 (2017), p.~e0176835.

\bibitem{logothetis2008we}
{\sc N.~K. Logothetis}, {\em What we can do and what we cannot do with {fMRI}}, Nature, 453 (2008), pp.~869--878.

\bibitem{menon202320}
{\sc V.~Menon}, {\em 20 years of the default mode network: A review and synthesis}, Neuron, 111 (2023), pp.~2469--2487.

\bibitem{niso2016omega}
{\sc G.~Niso, C.~Rogers, J.~T. Moreau, L.-Y. Chen, C.~Madjar, S.~Das, E.~Bock, F.~Tadel, A.~C. Evans, P.~Jolicoeur, et~al.}, {\em {OMEGA}: the open {MEG} archive}, Neuroimage, 124 (2016), pp.~1182--1187.

\bibitem{pascarella2019inversion}
{\sc A.~Pascarella, F.~Pitolli, et~al.}, {\em An inversion method based on random sampling for real-time meg neuroimaging}, Communications in Applied and Industrial Mathematics, 10 (2019), pp.~25--34.

\bibitem{raichle2009brief}
{\sc M.~E. Raichle}, {\em A brief history of human brain mapping}, Trends in neurosciences, 32 (2009), pp.~118--126.

\bibitem{raichle2015brain}
\leavevmode\vrule height 2pt depth -1.6pt width 23pt, {\em The brain's default mode network}, Annual review of neuroscience, 38 (2015), pp.~433--447.

\bibitem{salehi2020there}
{\sc M.~Salehi, A.~S. Greene, A.~Karbasi, X.~Shen, D.~Scheinost, and R.~T. Constable}, {\em There is no single functional atlas even for a single individual: Functional parcel definitions change with task}, NeuroImage, 208 (2020), p.~116366.

\bibitem{tadel2011brainstorm}
{\sc F.~Tadel, S.~Baillet, J.~C. Mosher, D.~Pantazis, and R.~M. Leahy}, {\em Brainstorm: A user-friendly application for meg/eeg analysis}, Computational intelligence and neuroscience, 2011 (2011), p.~879716.

\bibitem{uutela1999visualization}
{\sc K.~Uutela, M.~H{\"a}m{\"a}l{\"a}inen, and E.~Somersalo}, {\em Visualization of magnetoencephalographic data using minimum current estimates}, NeuroImage, 10 (1999), pp.~173--180.

\bibitem{wolters2004efficient}
{\sc C.~H. Wolters, L.~Grasedyck, and W.~Hackbusch}, {\em Efficient computation of lead field bases and influence matrix for the fem-based eeg and meg inverse problem}, Inverse problems, 20 (2004), p.~1099.

\bibitem{yuan2006model}
{\sc M.~Yuan and Y.~Lin}, {\em Model selection and estimation in regression with grouped variables}, Journal of the Royal Statistical Society Series B: Statistical Methodology, 68 (2006), pp.~49--67.

\end{thebibliography}
\end{document}